# PENALIZED MAXIMUM LIKELIHOOD AND SEMIPARAMETRIC SECOND-ORDER EFFICIENCY

By A. S. Dalalyan, G. K. Golubev and A. B. Tsybakov

*Université Paris VI, Université Paris VI and Université Aix-Marseille 1*

We consider the problem of estimation of a shift parameter of an unknown symmetric function in Gaussian white noise. We introduce a notion of semiparametric second-order efficiency and propose estimators that are semiparametrically efficient and second-order efficient in our model. These estimators are of a penalized maximum likelihood type with an appropriately chosen penalty. We argue that second-order efficiency is crucial in semiparametric problems since only the second-order terms in asymptotic expansion for the risk account for the behavior of the "nonparametric component" of a semiparametric procedure, and they are not dramatically smaller than the first-order terms.

**1. Introduction.** Semiparametric statistical models are the ones containing a finite-dimensional parameter of interest $\theta$ and an infinite-dimensional nuisance parameter $f$ which is a member of some large functional class. The goal is then to estimate $\theta$ efficiently without knowing $f$. A comprehensive account of the theory of semiparametric estimation is given in the book of Bickel, Klaassen, Ritov and Wellner [3]. In particular, it is shown that for many semiparametric models there exist estimators attaining the same asymptotic performance as efficient parametric estimators constructed for the problem where $f$ is completely specified. In other words, for such semiparametric models there is no loss of efficiency as compared to the corresponding parametric models with known $f$. These semiparametric models are usually called adaptive, but we prefer here to call them *S-adaptive*, or semiparametrically adaptive, in order to avoid confusion with nonparametric adaptivity to unknown smoothness of $f$. Estimators attaining parametric efficiency in $S$-adaptive models will be called $S$-adaptive (or efficient) estimators. Here and in what follows efficiency is understood in a local asymptotic minimax sense.

There exist various methods of constructing $S$-adaptive estimators. A general feature of these methods is that they proceed by "eliminating" the nonparametric component $f$, thus reducing the original semiparametric problem to a suitably chosen parametric one. The most common approach is to specify a least favorable parametric submodel of the full semiparametric model, locally in a neighborhood of $f$, and to estimate $\theta$ in such a submodel ([3, 22, 24, 30, 31, 32] and the references cited therein). Least favorable parametric submodels turn out to depend on $f$ only via a score function. "Elimination" of $f$ under this approach means to estimate nonparametrically the efficient score function. Resulting estimators of $\theta$ are often defined via one-step procedures that involve preliminary estimators of $\theta$ and nonparametric estimators of the efficient score function. We note here, in connection with the discussion that follows below, that results on efficiency and $S$-adaptivity are not very sensitive to the choice of preliminary nonparametric estimates of the efficient score function. For example, kernel, orthogonal series, nonparametric maximum likelihood and other estimates can be used, under rather wide assumptions on their parameters, such as kernels, bandwidths, etc. The important question of how to choose these parameters in practice is left open. Among other approaches that allow one to eliminate $f$ efficiently we mention profile likelihood techniques [25] and invariance-based inference [13].

Thus, for a variety of semiparametric models, the statistician actually has an entire library of efficient ($S$-adaptive) estimators of $\theta$. Which estimator is the best one? The theory discussed above does not answer this question because it deals only with the first-order asymptotics, which is the same for all $S$-adaptive estimators in a given model. Distinguishing between these estimators is possible on the basis of higher-order asymptotics. This motivates us to study here second-order asymptotics and second-order semiparametric efficiency. We would like to emphasize that a study of second-order effects is more important for semiparametric models than for purely parametric ones and it is crucial for practical implementation, at least for the following reasons.

- This is a compelling way to distinguish between various efficient semiparametric methods and to choose the best among them. More specifically, it allows one to choose optimally the smoothing parameters that define the "nonparametric component" of a given family of efficient semiparametric procedures.
- Second-order terms in asymptotics for semiparametric estimators are not dramatically smaller than the first-order terms; they might be in fact quite comparable to each other for moderate sample sizes. Second-order terms depend on the smoothness of $f$. For example, in a typical case of twice differentiable $f$ we get second-order terms $\sim n^{-7/10}$, the first-order asymptotics being as usual $n^{-1/2}$, where $n$ is the sample size. This differs from the purely parametric situation where the second-order terms decrease as $n^{-1}$ (cf. [20]).

Whereas first-order efficiency considerations for semiparametric models are essentially of a parametric flavor, the second-order ones come from nonparametric function estimation. Therefore, it is not surprising that the importance of second-order semiparametric asymptotics was first realized in the literature on nonparametric smoothing. Härdle and Tsybakov [15] pointed out that, in the single index model, the second-order term of the risk of the average derivative estimator is not significantly smaller than the first-order one and suggested choosing the optimal bandwidth by minimizing an asymptotic approximation of the second-order term. Mammen and Park [21] proceeded in a similar way to derive the optimal bandwidth for estimation of the efficient score function in the symmetric location problem. These papers considered specific families of estimators and did not deal with second-order efficiency among all estimators. Golubev and Härdle [9, 10] studied partial linear models and suggested second-order efficient estimators as well as their nonparametrically adaptive versions. These results rely strongly on the linearity and additivity of the parametric component in partial linear models. The problem of how to treat second-order efficiency for essentially nonlinear models has remained open, and our aim here is to give a solution to this problem.

We restrict our study to one basic model that seems to capture the main difficulties in deriving second-order efficiency, being at the same time simple enough to avoid unnecessary technicalities. Namely, we consider the estimation of a shift parameter $\theta$ based on observations

$$(1) \qquad x^{\varepsilon}(t) = f(t - \theta) + \varepsilon n(t), \qquad t \in [-1/2, 1/2],$$

where $n(t)$ is the standard Gaussian white noise process on $[-1/2, 1/2]$ (cf. [16], Chapter 3) and $f(\cdot)$ is a smooth symmetric [i.e., $f(t) = f(-t), \forall t$] periodic function with period 1, and $0 < \varepsilon < 1$ is a known noise parameter. With $\varepsilon = 1/\sqrt{n}$, where $n$ is an equivalent sample size, model (1) can be viewed as a "Gaussian white noise analog" of the classical symmetric location model [2, 26, 27].

If the signal $f$ is known, the maximum likelihood estimator

$$\hat{\theta}_{\mathrm{ML}} = \arg\max_{\tau} \int_{-1/2}^{1/2} f(t - \tau) x^{\varepsilon}(t) \, dt$$

is locally asymptotically minimax (e.g., [17]). In particular, its mean square risk satisfies

(2) $$\lim_{\varepsilon \to 0} \sup_{\theta \in \Theta} \mathbf{E}_{\theta, f}[(\hat{\theta}_{\mathrm{ML}} - \theta)^2 I^{\varepsilon}(f)] = 1,$$

for any sufficiently small interval $\Theta$, where

$$I^{\varepsilon}(f) = \varepsilon^{-2} \int_{-1/2}^{1/2} [f'(t)]^2 \, dt$$

is the Fisher information associated with model (1) and $\mathbf{E}_{\theta,f}$ is the expectation with respect to the distribution of the observation $\mathbf{X}^{\varepsilon} = \{x^{\varepsilon}(t), t \in [-1/2, 1/2]\}$ in model (1). The corresponding probability measure will be denoted by $\mathbf{P}_{\theta,f}$.

In a semiparametric setup where $f$ is not known, an efficient and $S$-adaptive estimator of $\theta$ is suggested by Golubev [8] for a model close to (1) where the observations are available for all $t \in \mathbf{R}$ and $f$ is not periodic. Härdle and Marron [14] discussed semiparametric estimation for models with discrete observations similar to (1) involving also a scale parameter.

Here we construct an $S$-adaptive and second-order efficient semiparametric estimator of $\theta$ in model (1). It is of penalized maximum likelihood type with an appropriately chosen penalty. To derive this estimator, we introduce a prior on $f$ and then maximize both in $\theta$ and $f$ the posterior density of $f$ given the observations. This procedure is of a Bayesian type w.r.t. $f$ for fixed $\theta$. It can be viewed in the following way: we "eliminate" the nonparametric component using a Bayesian argument, while the final estimation of $\theta$ is realized by maximum likelihood.

We conjecture that the penalized maximum likelihood approach using similar arguments would be a proper tool to get second-order efficient estimators for other semiparametric models, and we believe that our technique of proving minimax lower bounds with second-order terms might be useful there as well.

This paper is organized as follows. In Section 2 we give some heuristics concerning the first- and second-order efficiency in model (1). Section 3 contains the argument leading to a class of estimators defined by a sequence of weights: we show how these estimators (that are of penalized maximum likelihood type) are derived from Bayesian considerations. In Section 4 we show that, under certain assumptions on the sequence of weights, the estimators from this class are $S$-adaptive and we study their second-order asymptotics. Section 5 discusses a minimax problem for the second-order term. In Section 6 we give a locally asymptotically minimax lower bound and suggest a second-order efficient estimator obtained with a particular choice of weights. Sections 7–9 contain the proofs.

**2. Some heuristics.** This section provides some useful heuristics about first- and second-order semiparametric efficiency in model (1).

We first explain the result (2) obtained for known $f$. An intuitive way to do this is based on a local linear approximation of the signal $f(t-\theta)$. Suppose that $\theta$ belongs to a small interval $[\theta_0 - \Delta_\varepsilon, \theta_0 + \Delta_\varepsilon]$, where $\Delta_\varepsilon > 0$ and $\theta_0$ are known and $\Delta_\varepsilon \to 0$ as $\varepsilon \to 0$. This assumption is essentially equivalent to the existence of a $\Delta_\varepsilon$-consistent estimator of $\theta$. For simplicity, we assume that $\Delta_\varepsilon \sim \varepsilon$ [for rigorous proofs one needs to take $\Delta_\varepsilon$ slightly larger than $\varepsilon$, so that $\Delta_\varepsilon/\varepsilon \to \infty$, as $\varepsilon \to 0$, e.g., $\Delta_\varepsilon = \varepsilon\sqrt{\log(\varepsilon^{-2})}$]. Then, replacing $f(t-\theta)$ in (1) by its linear approximation $f(t-\theta_0) - f'(t-\theta_0)(\theta-\theta_0)$, we get the linear model

$$(3) \quad x_L^\varepsilon(t) = f(t-\theta_0) - f'(t-\theta_0)(\theta-\theta_0) + \varepsilon n(t), \qquad t \in [-1/2, 1/2].$$

When $f$ is known we can subtract $f(t-\theta_0)$ from these observations, thus obtaining an equivalent model,

$$y^\varepsilon(t) = f'(t-\theta_0)(\theta-\theta_0) + \varepsilon n(t), \qquad t \in [-1/2, 1/2].$$

Estimation of $\theta - \theta_0$ in this linear regression model is now straightforward. Multiplying the observation $y^\varepsilon(t)$ by $f'(t-\theta_0)$, integrating over the interval $[-1/2, 1/2]$ and dividing by $I^\varepsilon(f)$ we get the Gaussian shift model

$$(4) \qquad Y^\varepsilon = \theta - \theta_0 + [I^\varepsilon(f)]^{-1/2}\xi,$$

where $\xi \sim \mathcal{N}(0,1)$. Clearly, $Y^\varepsilon$ is an efficient estimator of $\theta - \theta_0$. Thus, the argument here is based on replacing the original nonlinear estimation problem by a Gaussian shift experiment. A deep theoretical background for this argument is given by Le Cam's theory of asymptotic equivalence [19].

Suppose now that $f$ is an unknown symmetric function. Then again we can use model (3) to approximate the initial model (1). But the approximating model is now nonlinear since it contains the product of unknown parameters $(\theta - \theta_0)$ and $f'(t-\theta_0)$. Fortunately, this is not a problem, and in this case one can also construct an efficient estimator.

Indeed, since $f'$ is an odd function and $f$ is an even function, projecting the observations (3) on the spaces of even and odd functions we get

$$(5) \qquad x_e^\varepsilon(t) = f(t-\theta_0) + \varepsilon n_e(t),$$

$$(6) \qquad x_o^\varepsilon(t) = f'(t-\theta_0)(\theta-\theta_0) + \varepsilon n_o(t),$$

where $n_o(t)$ and $n_e(t)$ are two independent Gaussian white noise processes. Based on $x_e^\varepsilon(t)$, we estimate the derivative $f'(t-\theta_0)$ and then plug this estimator into (6) to recover the parameter of interest from the observation $x_o^\varepsilon(t)$. This allows us to obtain an efficient ($S$-adaptive) estimator of $\theta$.

We turn now to a heuristic derivation of second-order asymptotics. In order to do that we simplify our approximate statistical model (5)–(6) assuming that $\theta_0 = 0$ and translating the observations $x_{\rm o}^\varepsilon(t), x_{\rm e}^\varepsilon(t)$ in a sequence space.

We will suppose throughout the paper that the unknown function $f$ can be represented as

$$f(t) = \sqrt{2} \sum_{k=1}^{\infty} f_k \cos(2\pi k t), \tag{7}$$

where the Fourier series converges for all $t$ and the Fourier coefficients $f_k$ are defined by

$$f_k = \sqrt{2} \int_{-1/2}^{1/2} f(t) \cos(2\pi k t)\, dt.$$

Using this and projecting (5) and (6) on the trigonometric basis functions we obtain the sequence model

$$X_k = f_k + \varepsilon \xi_k, \qquad k = 1, 2, \ldots, \tag{8}$$

$$X_k^* = \theta(2\pi k) f_k + \varepsilon \xi_k^*, \qquad k = 1, 2, \ldots, \tag{9}$$

where $(\xi_k, \xi_k^*,\ k = 1, 2, \ldots)$ are i.i.d. standard Gaussian random variables. The nuisance parameters $f_k$ can be estimated from (8) by well-known techniques for the Gaussian sequence model (see, e.g., [29]). In particular, it is natural to use linear estimators of $f_k$ defined by $\hat{f}_k = h_k X_k$, where $h_k = h_k(\varepsilon)$ are such that $\sum_{k=1}^{\infty} h_k^2 < \infty$. An example is $h_k = \mathbb{1}_{\{k \le N_\varepsilon\}}$ where $\mathbb{1}_{\{\cdot\}}$ is the indicator function and $N_\varepsilon$ is an integer such that $N_\varepsilon \to \infty$ as $\varepsilon \to 0$.

Next, considering separately model (9), it is not hard to show that if $f_k$ were known the maximum likelihood (least squares) estimator

$$\hat{\theta}_f = \sum_{k=1}^{\infty} (2\pi k) f_k X_k^* \Big/ \sum_{k=1}^{\infty} (2\pi k)^2 f_k^2 \tag{10}$$

would be asymptotically minimax for $\theta$. At first sight, it seems natural to plug in $\hat{f}_k$ instead of $f_k$ in the expression for $\hat{\theta}_f$, thus obtaining the estimator

$$\tilde{\theta} = \sum_{k=1}^{\infty} (2\pi k) h_k X_k X_k^* \Big/ \sum_{k=1}^{\infty} (2\pi k)^2 h_k^2 X_k^2. \tag{11}$$

However, this estimator is not optimal: it can have a very large bias. The reason is that the functional $\sum_{k=1}^{\infty} (2\pi k)^2 f_k^2$ in (10) is not estimated correctly. An improved version of $\tilde{\theta}$ can be suggested in the form

$$\theta^* = \sum_{k=1}^{\infty} (2\pi k) h_k X_k X_k^* \Big/ \sum_{k=1}^{\infty} (2\pi k)^2 h_k (X_k^2 - \varepsilon^2). \tag{12}$$

As compared to (11), we replace $h_k^2$ by $h_k$ in the denominator and replace $X_k^2$ by the unbiased estimator $X_k^2 - \varepsilon^2$ of $f_k^2$. This turns out to improve significantly the asymptotics of the risk.

We now give a heuristic analysis of the risk of $\theta^*$. Using (8)–(9) and the notation $\|f'\|^2 = \varepsilon^2 I^\varepsilon(f) = \sum_{k=1}^\infty (2\pi k)^2 f_k^2$, we obtain

$$(13) \qquad (\theta^* - \theta)\sqrt{I^\varepsilon(f)} = \|f'\| \frac{\chi^\varepsilon - \Gamma_1^\varepsilon}{\sum_{k=1}^\infty (2\pi k)^2 h_k f_k^2 + \Gamma_2^\varepsilon},$$

where

$$\chi^\varepsilon = \sum_{k=1}^\infty (2\pi k) h_k f_k \xi_k^* + \varepsilon \sum_{k=1}^\infty (2\pi k) h_k \xi_k \xi_k^*,$$

$$\Gamma_1^\varepsilon = \theta \sum_{k=1}^\infty (2\pi k)^2 h_k f_k \xi_k + \theta \varepsilon \sum_{k=1}^\infty (2\pi k)^2 h_k (\xi_k^2 - 1),$$

$$\Gamma_2^\varepsilon = 2\varepsilon \sum_{k=1}^\infty (2\pi k)^2 h_k f_k \xi_k + \varepsilon^2 \sum_{k=1}^\infty (2\pi k)^2 h_k (\xi_k^2 - 1).$$

In order to simplify the expression in (13) we assume that $\sum_{k=1}^\infty (2\pi k)^4 f_k^2 < \infty$ and that $h_k$ are chosen so that $\varepsilon \sum_{k=1}^\infty (2\pi k)^2 h_k < \infty$. Under these conditions, using $|\theta| \leq \Delta_\varepsilon \sim \varepsilon$, one obtains that

$$\mathbf{E}_{\theta,f}[(\Gamma_1^\varepsilon)^2] = O(\varepsilon^2), \qquad \mathbf{E}_{\theta,f}[(\Gamma_2^\varepsilon)^2] = O(\varepsilon^2).$$

It is also straightforward to see that $\mathbf{E}_{\theta,f}(\chi^\varepsilon \Gamma_1^\varepsilon) = 0$, and to show, with some easy algebra, that

$$\mathbf{E}_{\theta,f}[(\chi^\varepsilon)^2 \Gamma_2^\varepsilon] = 4\varepsilon^3 \sum_{k=1}^\infty h_k^3 (2\pi k)^4 f_k^2 + 2\varepsilon^4 \sum_{k=1}^\infty h_k^3 (2\pi k)^4 = O(\varepsilon^2).$$

Next note that we are allowed to drop the terms of order $O(\varepsilon^2)$ since their contribution in the risk (asymptotically, in the mean absolute value) is smaller than the final second-order asymptotics that we are going to obtain. Up to these terms, we get from (13)

$$(\theta^* - \theta)\sqrt{I^\varepsilon(f)} \approx \frac{\|f'\|}{\sum_{k=1}^\infty (2\pi k)^2 h_k f_k^2} \times \left[\chi^\varepsilon - \Gamma_1^\varepsilon - \chi^\varepsilon \Gamma_2^\varepsilon \left(\sum_{k=1}^\infty (2\pi k)^2 h_k f_k^2\right)^{-1}\right],$$

and thus

$$\mathbf{E}_{\theta,f}[(\theta^* - \theta)^2 I^\varepsilon(f)] \approx \|f'\|^2 \left(\sum_{k=1}^\infty (2\pi k)^2 h_k f_k^2\right)^{-2} \mathbf{E}_{\theta,f}[(\chi^\varepsilon)^2]$$

$$= \|f'\|^2 \sum_{k=1}^\infty (2\pi k)^2 h_k^2 (\varepsilon^2 + f_k^2) \left(\sum_{k=1}^\infty (2\pi k)^2 h_k f_k^2\right)^{-2}.$$

This expression can be simplified if we assume that $0 \leq h_k \leq 1$ and

$$\left[\sum_{k=1}^{\infty}(1-h_k)(2\pi k)^2 f_k^2\right]^2 = o\left(\sum_{k=1}^{\infty}(1-h_k)^2(2\pi k)^2 f_k^2\right) \tag{14}$$

as $\varepsilon \to 0$. Then, in particular, $\sum_{k=1}^{\infty}(1-h_k)(2\pi k)^2 f_k^2 = o(1)$, and one obtains

$$\left(\sum_{k=1}^{\infty}(2\pi k)^2 h_k f_k^2\right)^{-2} = \left(\|f'\|^2 + \sum_{k=1}^{\infty}(2\pi k)^2(h_k-1)f_k^2\right)^{-2}$$

$$\approx \|f'\|^{-4}\left[1 - 2\|f'\|^{-2}\sum_{k=1}^{\infty}(2\pi k)^2(h_k-1)f_k^2\right].$$

Using this and (14) we derive the following expansion for the risk:

$$\mathbf{E}_{\theta,f}[(\theta^*-\theta)^2 I^\varepsilon(f)] \approx \left[1 + \|f'\|^{-2}\left(\sum_{k=1}^{\infty}(2\pi k)^2(\varepsilon^2 h_k^2 + (h_k^2-1)f_k^2)\right)\right]$$

$$\times \left[1 - 2\|f'\|^{-2}\sum_{k=1}^{\infty}(2\pi k)^2(h_k-1)f_k^2\right] \tag{15}$$

$$\approx 1 + \|f'\|^{-2} R^\varepsilon[f,h],$$

where

$$R^\varepsilon[f,h] = \sum_{k=1}^{\infty}(2\pi k)^2[(1-h_k)^2 f_k^2 + \varepsilon^2 h_k^2]. \tag{16}$$

The second-order term in (15), that is, the functional $\|f'\|^{-2}R^\varepsilon[f,h]$, has a clear statistical meaning. Suppose that we know $\theta$ and we want to estimate the derivative $f'(t-\theta)$ based on observations (1). To measure the quality of an estimator $\hat{f}'(t-\theta)$ we choose the relative mean integrated squared error,

$$\mathrm{Err}(\hat{f}', f') = \frac{1}{\|f'\|^2}\mathbf{E}_{\theta,f}\int_{-1/2}^{1/2}[\hat{f}'(t-\theta) - f'(t-\theta)]^2\,dt.$$

Consider a linear estimator

$$\tilde{f}'_h(t-\theta) = -2\sum_{k=1}^{\infty} h_k(2\pi k)\sin[2\pi k(t-\theta)]\int_{-1/2}^{1/2}\cos[2\pi k(t-\theta)]x^\varepsilon(t)\,dt.$$

Using (7), it is easy to show that $\mathrm{Err}(\tilde{f}'_h, f') = \|f'\|^{-2}R^\varepsilon[f,h]$. Thus, the expression $\|f'\|^{-2}R^\varepsilon[f,h]$ is a relative mean integrated squared error for nonparametric estimation of the derivative of $f$ in the Gaussian white noise model. We see that the second-order expansion (15) relates two statistical problems: semiparametric estimation of $\theta$ and nonparametric estimation in $L_2$-norm of the Fisher informant $f'(t - \theta)$. It also reveals a presumably general fact that second-order asymptotic terms in semiparametric problems account for the mean integrated squared error of recovering of the Fisher informant.

**3. Penalized maximum likelihood estimator.** In Section 2 we have sketched second-order asymptotics for the estimator $\theta^*$ in model (8)–(9), which is only a local approximation of the original model (1) in a neighborhood of $\theta_0 = 0$. Thus, $\theta^*$ is not directly applicable for model (1). Of course, the procedure can be corrected: instead of replacing $\theta_0$ by 0, one should replace it by a preliminary $\varepsilon$-consistent estimator of $\theta$. This would lead to a two-stage estimation procedure that would presumably have the desired second-order behavior under some conditions. There exists, however, a direct and more elegant estimator achieving the same result. This estimator is inspired by the Bayes argument that we are going to describe now.

Given model (1), we have at our disposition the following series of discrete observations:

$$
\begin{aligned}
x_k &= f_k \cos(2\pi k\theta) + \varepsilon \xi_k, \\
x_k^* &= f_k \sin(2\pi k\theta) + \varepsilon \xi_k^*, \qquad k = 1, 2, \ldots.
\end{aligned}
\tag{17}
$$

Here $(\xi_k, \xi_k^*,\ k = 1, 2, \ldots)$ are i.i.d. standard Gaussian random variables,

$$
x_k = \sqrt{2} \int_{-1/2}^{1/2} x^\varepsilon(t) \cos(2\pi kt)\, dt, \qquad x_k^* = \sqrt{2} \int_{-1/2}^{1/2} x^\varepsilon(t) \sin(2\pi kt)\, dt,
$$

and (17) is obtained by projection of (1) on the trigonometric basis functions on $[-1/2, 1/2]$ using (7).

Our aim is to define a suitable estimator of $\theta$ using these observations. A general idea is to "eliminate" first the nonparametric component of the model represented by the sequence of Fourier coefficients $f_k$ (which we consider to be nuisance parameters). We will proceed as follows. Assume for a moment that the $f_k$'s are independent zero-mean Gaussian random variables with variances $\sigma_k^2$. Assume also that they are independent of the noise sequence $\{\xi_k, \xi_k^*\}$. We will replace the sequence $\{f_k\}$ by the most probable, with respect to the posterior distribution of $\{f_k\}$ given $\{x_k, x_k^*\}$, sequence $\{f_k^*\}$. Clearly, this sequence will depend only on $\{x_k, x_k^*\}$ and $\theta$, and thus $\{f_k\}$ will be eliminated. The final step will be to maximize over $\theta$ the remaining likelihood, thus obtaining an estimator of $\theta$.

To define the procedure formally, note that the problem factorizes: it is sufficient to find $f_k^*$'s for a fixed $k$, since the triples $x_k, x_k^*, f_k$ with different $k$ are independent. Maximizing over $f_k$ the posterior density of $f_k$ given $x_k, x_k^*$ is equivalent to maximizing the joint density of $x_k, x_k^*, f_k$, which equals

$$p_\theta(x_k, x_k^*, f_k) = \left(\frac{1}{\sqrt{2\pi}}\right)^3 \sigma_k^{-1} \exp\left[-\frac{f_k^2}{2\sigma_k^2}\right]$$

$$\times \exp\left[-\frac{(x_k - f_k \cos(2\pi k\theta))^2 + (x_k^* - f_k \sin(2\pi k\theta))^2}{2\varepsilon^2}\right]$$

$$= A(x_k, x_k^*)$$

$$\times \exp\left[\frac{\sqrt{2}f_k}{\varepsilon^2}\int_{-1/2}^{1/2} \cos[2\pi k(t-\theta)]x^\varepsilon(t)\,dt - \frac{f_k^2(\varepsilon^2 + \sigma_k^2)}{2\varepsilon^2\sigma_k^2}\right],$$

where $A(x_k, x_k^*)$ does not depend on $f_k$ and $\theta$. The maximizer of $p_\theta(x_k, x_k^*, f_k)$ over $f_k$ has the form

$$f_k^*(\theta) = \sqrt{2}\lambda_k \int_{-1/2}^{1/2} \cos[2\pi k(t-\theta)]x^\varepsilon(t)\,dt,$$

where $\lambda_k = \dfrac{\sigma_k^2}{\sigma_k^2 + \varepsilon^2}$ and

(18)
$$\max_{f_k} p_\theta(x_k, x_k^*, f_k)$$
$$= p_\theta(x_k, x_k^*, f_k^*(\theta))$$
$$= A(x_k, x_k^*)\exp\left[\lambda_k\left(\int_{-1/2}^{1/2}\cos[2\pi k(t-\theta)]x^\varepsilon(t)\,dt\right)^2\right].$$

Set

$$\hat\theta_{\text{PML}} = \arg\max_{\theta\in\Theta}\prod_{k=1}^\infty p_\theta(x_k, x_k^*, f_k^*(\theta)) = \arg\max_{\theta\in\Theta}\left[\max_{\{f_k\}}\prod_{k=1}^\infty p_\theta(x_k, x_k^*, f_k)\right],$$

where $\Theta$ is a parameter set associated with the model. Thus, $\hat\theta_{\text{PML}}$ is the $\theta$-component of the overall maximum likelihood estimator corresponding to the infinite product density $\prod_{k=1}^\infty p_\theta(x_k, x_k^*, f_k)$. In view of (18), we may write this estimator as

(19) $$\hat\theta_{\text{PML}} = \arg\max_{\tau\in\Theta}\left\{\sum_{k=1}^\infty \lambda_k\left(\int_{-1/2}^{1/2}\cos[2\pi k(t-\tau)]x^\varepsilon(t)\,dt\right)^2\right\},$$

or as

(20) $$\hat\theta_{\text{PML}} = \arg\max_{\tau\in\Theta}\max_{\{g_k\}}\left[\frac{\sqrt{2}}{\varepsilon^2}\sum_{k=1}^\infty g_k \int_{-1/2}^{1/2}\cos[2\pi k(t-\tau)]x^\varepsilon(t)\,dt\right.$$
$$\left.-\sum_{k=1}^\infty g_k^2\left(\frac{1}{2\varepsilon^2} + \frac{1}{2\sigma_k^2}\right)\right],$$

where $\max_{\{g_k\}}$ denotes the maximum over sequences $\{g_k\}$ belonging to a subset of $\ell_2$, and we suppose that $f$ satisfies conditions such that the infinite sums converge almost surely. We will call $\hat{\theta}_{\text{PML}}$ a penalized maximum likelihood estimator (PMLE), although this is not a PMLE in the usual sense. Comparing $\hat{\theta}_{\text{PML}}$ with the maximum likelihood estimator $\hat{\theta}_{\text{ML}}$, we see that $\hat{\theta}_{\text{PML}}$ can be interpreted as a penalized version of $\hat{\theta}_{\text{ML}}$ corresponding to a function $f(\cdot) = f_\tau(\cdot)$ whose Fourier coefficients are the maximizers $\{g_k^*(\tau)\}$ of the term in square brackets in (20) over $\{g_k\}$ for fixed $\tau$ and to the penalty $\sum_{k=1}^{\infty}(g_k^*(\tau))^2(\frac{1}{2\varepsilon^2} + \frac{1}{2\sigma_k^2})$ (up to a multiplicative constant, cf. definition of $\hat{\theta}_{\text{ML}}$). Thus, the difference of $\hat{\theta}_{\text{PML}}$ from the "pure" PMLE is in the fact that $f(\cdot) = f_\tau(\cdot)$ is not fixed and known: it depends on the parameter $\tau$ over which the maximization is carried out.

To make the estimator $\hat{\theta}_{\text{PML}}$ feasible, it is natural to consider only finite sums in (19), (20), including the terms with $k \leq N_\varepsilon$, for some $N_\varepsilon$ that depends on $\varepsilon$ and tends to $\infty$ as $\varepsilon \to 0$. In particular, this will be the case for the second-order minimax estimator that we derive below.

Note that the estimator (12) defined in Section 2 is nothing but a local version of the estimator (19) in a neighborhood of $\theta_0 = 0$. In fact, differentiating formally the expression in curly brackets in (19) we obtain that $\hat{\theta}_{\text{PML}}$ is a solution of the equation

$$
\begin{aligned}
(21) \quad \sum_{k=1}^{\infty} \lambda_k (2\pi k) &\left( \int \cos[2\pi k(t-\tau)] x^\varepsilon(t)\, dt \right) \\
&\times \left( \int \sin[2\pi k(t-\tau)] x^\varepsilon(t)\, dt \right) = 0.
\end{aligned}
$$

The integrals in (21) are equal to $y_k = x_k \cos(2\pi k\tau) + x_k^* \sin(2\pi k\tau)$ and $y_k' = x_k^* \cos(2\pi k\tau) - x_k \sin(2\pi k\tau)$, respectively, allowing one to reduce (21) to

$$\sum_{k=1}^{\infty} \lambda_k(2\pi k)\{x_k x_k^* \cos(4\pi k\tau) - [x_k^2 - (x_k^*)^2]\sin(4\pi k\tau)/2\} = 0.$$

Linearizing this equation in the vicinity of $\tau = 0$, we get the following approximate formula for a solution of (21):

$$(22) \quad \hat{\theta}_{\text{PML}} \approx \sum_{k=1}^{\infty}(2\pi k)\lambda_k x_k x_k^* \Big/ \sum_{k=1}^{\infty}(2\pi k)^2 \lambda_k [x_k^2 - (x_k^*)^2].$$

It can be shown, using the argument from Section 2 that (22) is asymptotically analogous to the estimator $\theta^*$ given by (12) with $h_k = \lambda_k$. One difference is that here we have $x_k, x_k^*$ instead of $X_k, X_k^*$, but $x_k \approx X_k$ and $x_k^* \approx X_k^*$ for $\theta$ close to 0. Another point is that these estimators have somewhat different denominators. However, for small $\theta$ both denominators estimate the same quadratic functional $\sum_{k=1}^{\infty}(2\pi k)^2 f_k^2$ and one can show that they are quite close to each other, so that their difference does not appear in the second-order asymptotics of the risk.

**4. Second-order asymptotics of the estimators.** In this section we consider the class of estimators defined by

$$(23) \qquad \hat{\theta}_{\mathrm{AD}} = \arg\max_{\tau \in \Theta} \left\{ \sum_{k=1}^{\infty} h_k \left( \int_{-1/2}^{1/2} \cos[2\pi k(t-\tau)] x^{\varepsilon}(t) \, dt \right)^2 \right\},$$

where $\{h_k\}$ is a sequence of real numbers satisfying some general conditions. For a particular choice $h_k = \lambda_k$ the estimator $\hat{\theta}_{\mathrm{AD}}$ is equal to the penalized maximum likelihood estimator (19) obtained from a Bayesian argument with $\lambda_k = \sigma_k^2/(\sigma_k^2 + \varepsilon^2)$, but we also allow other weights $h_k$. In particular, the weights $\{h_k\}$ such that $h_k = 1$ for some initial values of $k$ play an important role in our further argument, while we always have $\lambda_k < 1$ for $\hat{\theta}_{\mathrm{PML}}$.

We will show that under some assumptions on $\{h_k\}$ the estimator $\hat{\theta}_{\mathrm{AD}}$ is $S$-adaptive and we will give explicit second-order asymptotics for the risk of $\hat{\theta}_{\mathrm{AD}}$. In what follows we will suppose that $h_k \neq 0$ for only a finite (typically, depending on $\varepsilon$ and growing to $\infty$, as $\varepsilon \to 0$) number of integers $k$. This assumption is natural, since otherwise the estimator $\hat{\theta}_{\mathrm{AD}}$ is not feasible. In order not to specify the set where $h_k \neq 0$ we keep in the notation the sums over all integers $k$.

We first define the parametric set $\Theta$ where $\theta$ lies. Since $f$ is symmetric and periodic with period 1, we get that $s(t) = f(1/2 - t)$ is also symmetric and periodic with period 1. Hence, the observations $x^{\varepsilon}(t)$ corresponding to parameters $(\theta, f(\cdot))$ and $(\theta - 1/2, s(\cdot))$ have the same probability distribution. So we cannot discriminate between values $\theta, \theta + 1/2, \theta + 1, \ldots$ in model (1) if we suppose that $f$ belongs to the class of symmetric and periodic functions with period 1. In order that the model be identifiable, $\Theta$ should be strictly included in an interval of length $1/2$. For definiteness, we assume the following.

ASSUMPTION A1.   $\Theta = \{\theta : |\theta| \leq \tau_0\}$ where $0 < \tau_0 < 1/4$.

Next, we define the class of functions $F$ where $f$ lies. Let $\rho$ and $C_0$ be positive constants. Denote by $F = F(\rho, C_0)$ the class of all functions $f : [-1/2, 1/2] \to \mathbf{R}$ that admit the Fourier expansion (7) with coefficients $f_k$ satisfying the following assumptions.

ASSUMPTION A2.   $f_1^2 \geq \rho$.

ASSUMPTION A3.   $\|f''\|^2 \leq C_0$.

Here and in the sequel, for a sequence of real numbers $\{a_k\}$, we use the notation

$$\|a\|^2 = \sum_{k=1}^{\infty} a_k^2, \qquad \|a'\|^2 = \sum_{k=1}^{\infty} a_k^2 (2\pi k)^2, \qquad \|a''\|^2 = \sum_{k=1}^{\infty} a_k^2 (2\pi k)^4.$$

Assumptions A2 and A3 imply that

(24) $$C_0 \geq \|f'\|^2 \geq (2\pi)^2 \rho \qquad \forall f \in F.$$

Furthermore, we impose some conditions on the weight sequence $\{h_k\}$, assuming that it depends on $\varepsilon$.

ASSUMPTION B. The weight sequence $\{h_k\}$ is such that $h_1 = 1$, $0 \leq h_k \leq 1$ for all $k$, and

B1. $\|h'\| \geq \rho_1 \log^2(\varepsilon^{-2}) \max_k h_k(2\pi k)$, where $\rho_1 > 0$ is a constant that does not depend on $\varepsilon$,
B2. $\varepsilon^2 \sum_{k=1}^{\infty} h_k (2\pi k)^4 \leq C_1$, where $C_1$ is a constant that does not depend on $\varepsilon$.

We remark that the condition $0 \leq h_k \leq 1$ here is quite natural: if $h_k \notin [0,1]$, projecting $h_k$ on $[0,1]$ only improves the second-order asymptotics (cf. the expression for $R^\varepsilon[f,h]$ in (16)). Note also that Assumption B2 and the fact that $0 \leq h_k \leq 1$ imply the finiteness of $\|h'\|$ for any $\varepsilon$. Assumptions B1 and B2 are not very restrictive. For example, consider the projection weights $h_k = \mathbb{1}_{\{k \leq N_\varepsilon\}}$ where $N_\varepsilon$ is an integer such that $N_\varepsilon \to \infty$ as $\varepsilon \to 0$. Then Assumption B1 is equivalent to $N_\varepsilon \geq C \log^4(\varepsilon^{-2})$ for some constant $C > 0$, and Assumption B2 is satisfied if $N_\varepsilon = O(\varepsilon^{-2/5})$ as $\varepsilon \to 0$.

Finally, we will need the following assumption involving both $f$ and $\{h_k\}$.

ASSUMPTION C. The weight sequence $\{h_k\}$ is such that, uniformly in $f \in F$,

$$\left[\sum_{k=1}^{\infty} (1-h_k)(2\pi k)^2 f_k^2\right]^2 = o\left(\sum_{k=1}^{\infty} (1-h_k)^2 (2\pi k)^2 f_k^2\right) \qquad \text{as } \varepsilon \to 0.$$

Note that, again, Assumption C is quite mild. For the projection weights $h_k = \mathbb{1}_{\{k \leq N_\varepsilon\}}$ it means that $\sum_{k \geq N_\varepsilon} (2\pi k)^2 f_k^2 \to 0$ as $\varepsilon \to 0$, uniformly in $f \in F$, which is true due to Assumption A3.

THEOREM 1. *Let Assumptions* A1–A3, B *and* C *be satisfied. Then, uniformly in* $f \in F$ *and in* $\theta \in \Theta$,

$$\mathbf{E}_{\theta,f}[(\hat{\theta}_{\mathrm{AD}} - \theta)^2 I^\varepsilon(f)] = 1 + (1 + o(1)) \frac{R^\varepsilon[f,h]}{\|f'\|^2} \qquad \text{as } \varepsilon \to 0,$$

*where the functional* $R^\varepsilon[f,h]$ *is defined in* (16).

Proof of Theorem 1 is given in Section 7.

Assumptions A3, B and C imply that

$$\sup_{f \in F} R^\varepsilon[f, h] = o(1) \qquad \text{as } \varepsilon \to 0.$$

In fact, it follows from Assumptions B1 and B2 that $\varepsilon^2 \sum_{k=1}^\infty h_k^2 (2\pi k)^2 = o(1)$, while Assumptions A3 and C yield $\sum_{k=1}^\infty (1-h_k)^2 (2\pi k)^2 f_k^2 = o(1)$ as $\varepsilon \to 0$. Thus, Theorem 1 shows that $\hat{\theta}_{\text{AD}}$ has the same first-order asymptotics as the efficient estimator $\hat{\theta}_{\text{ML}}$ [cf. (2)], that is, $\hat{\theta}_{\text{AD}}$ is $S$-adaptive under the assumptions of Theorem 1. But Theorem 1 says more than that, because it also provides an asymptotically exact second-order expansion for the risk of $\hat{\theta}_{\text{AD}}$.

**5. Minimax problem for second-order term.** It follows from Theorem 1 that the second-order term of the risk of $\hat{\theta}_{\text{AD}}$ depends on the coefficients $\{h_k\}$ only via the functional $R^\varepsilon[f, h]$. We would like to make this term as small as possible by minimizing it over $h_k$. Since we do not know the nuisance parameters $f_k$ we consider a minimax setting: we look for $h = \{h_k\}$ that minimizes the maximum of the functional $R^\varepsilon[f, h]$ over a suitably chosen set of sequences $\{f_k\}$. Namely, we consider a Sobolev ball

$$\mathcal{W}(\beta, L) = \left\{ f : \sum_{k=1}^\infty (2\pi k)^{2\beta} f_k^2 \leq L \right\},$$

where $\beta > 1$ and $L > 0$ are given constants. A minimax sequence of weights $q = \{q_k\} \in \ell_2$ is defined by

$$\sup_{f \in \mathcal{W}(\beta, L)} R^\varepsilon[f, q] = \inf_{h \in \ell_2} \sup_{f \in \mathcal{W}(\beta, L)} R^\varepsilon[f, h].$$

It is well known (see, e.g., [1] or [23]) that such a sequence $q$ exists and it has the form

$$(25) \qquad q_k = \left[ 1 - \left( \frac{k}{W_\varepsilon} \right)^{\beta - 1} \right]_+,$$

where $x_+ = \max(x, 0)$ and $W_\varepsilon$ is a solution of the equation

$$(26) \qquad \varepsilon^2 \sum_{k=1}^\infty \left[ \left( \frac{W_\varepsilon}{k} \right)^{\beta - 1} - 1 \right]_+ (2\pi k)^{2\beta} = L.$$

As $\varepsilon \to 0$, we have

$$(27) \qquad W_\varepsilon = (1 + o(1)) \left( \frac{L}{\varepsilon^2} \frac{(\beta + 2)(2\beta + 1)}{(2\pi)^{2\beta}(\beta - 1)} \right)^{1/(2\beta + 1)}.$$

Moreover, the functional $R^\varepsilon[f, h]$ has a saddle point on $\mathcal{W}(\beta, L) \times \ell_2$ (cf. [1], [23] or [29], Chapter 3) with components $s, q$, where $s = \{s_k\}$ is any sequence satisfying

$$(28) \qquad s_k^2 = \varepsilon^2 \frac{q_k}{1 - q_k} = \varepsilon^2 \left[\left(\frac{W_\varepsilon}{k}\right)^{\beta-1} - 1\right]_+.$$

The existence of a saddle point at $(s, q)$ means that

$$\inf_{h \in \ell_2} \sup_{f \in \mathcal{W}(\beta, L)} R^\varepsilon[f, h] = \sup_{f \in \mathcal{W}(\beta, L)} \inf_{h \in \ell_2} R^\varepsilon[f, h] = R^\varepsilon[s, q].$$

Using (25), (26) and (28), the value $R^\varepsilon[s, q]$ can be expressed explicitly, which yields

$$(29) \quad \inf_{h \in \ell_2} \sup_{f \in \mathcal{W}(\beta, L)} R^\varepsilon[f, h] = \sup_{f \in \mathcal{W}(\beta, L)} R^\varepsilon[f, q] = \varepsilon^2 \sum_{k=1}^{\infty} (2\pi k)^2 q_k \stackrel{\text{def}}{=} r^\varepsilon.$$

Note finally that, as $\varepsilon \to 0$,

$$(30) \quad \begin{aligned} r^\varepsilon &= \frac{(2\pi)^2(\beta - 1)}{3(\beta + 2)} \varepsilon^2 W_\varepsilon^3 (1 + o(1)) \\ &= C^*(\beta, L) \varepsilon^{(4\beta-4)/(2\beta+1)} (1 + o(1)), \end{aligned}$$

where

$$C^*(\beta, L) = \frac{1}{3} \left(\frac{\beta - 1}{2\pi(\beta + 2)}\right)^{(2\beta-2)/(2\beta+1)} (L(2\beta + 1))^{3/(2\beta+1)}.$$

The rate $\varepsilon^{(4\beta-4)/(2\beta+1)}$ in (30) characterizes the ratio of second-order terms to first-order terms in the asymptotic expansion for the nonnormalized risk $\mathbf{E}_{\theta,f}[(\hat{\theta}_\varepsilon - \theta)^2]$. This ratio is not dramatically small for $\beta$ not too large; for example, it equals $\varepsilon^{4/5}$ for $\beta = 2$. Thus, the second-order terms might be comparable with the first-order ones. In absolute value, the first-order term of nonnormalized risk decreases as $\varepsilon^2$ and the second-order term as $\varepsilon^{(8\beta-2)/(2\beta+1)}$.

**6. Locally minimax lower bound and second-order efficiency.** In this section we obtain a lower bound for the minimax risk and construct a second-order efficient estimator of $\theta$.

Let $\bar{f}$ be a fixed function from $F(\rho, C_0)$ with the Fourier coefficients denoted by $\bar{f}_k$. For $\delta > 0$ define a vicinity of $\bar{f}$ by

$$(31) \qquad F_\delta(\bar{f}) = \{f = \bar{f} + v : \|v\| \leq \delta, v \in \mathcal{W}(\beta, L)\}.$$

It is assumed that $\beta \geq 2$. Recall that $\|\bar{f}''\| < \infty$ since $\bar{f} \in F(\rho, C_0)$ (cf. Assumption A3). If $\delta$ is small enough, $F_\delta(\bar{f}) \subseteq F(\rho', C_0')$ for some $\rho' > 0$, $C_0' > 0$ depending only on $\rho, C_0, L$.

THEOREM 2. *Let the real number $\delta = \delta_\varepsilon$ be such that $\lim_{\varepsilon \to 0} \delta_\varepsilon = 0$ and $\lim_{\varepsilon \to 0} \delta_\varepsilon^2 / (\varepsilon^2 W_\varepsilon^{1+\alpha}) = \infty$ for some $\alpha > 0$, where $W_\varepsilon$ satisfies (27). Then, as $\varepsilon \to 0$,*

$$(32) \qquad \inf_{\hat{\theta}_\varepsilon} \sup_{\theta \in \Theta, f \in F_{\delta_\varepsilon}(\bar{f})} \mathbf{E}_{\theta,f}[(\hat{\theta}_\varepsilon - \theta)^2 I^\varepsilon(f)] \geq 1 + (1 + o(1)) \frac{r^\varepsilon}{\|\bar{f}'\|^2}.$$

*Here and in what follows $\inf_{\hat{\theta}_\varepsilon}$ (or $\inf_{\hat{\theta}}$) is the infimum over all estimators based on the observation $\mathbf{X}^\varepsilon$, and $r^\varepsilon$ is the minimax value defined in* (29).

The proof of Theorem 2 is given in Section 8.

Motivated by the above results, we introduce the following notion of semiparametric second-order efficiency.

DEFINITION 1. An estimator $\theta_\varepsilon^*$ is called second-order efficient at $\bar{f} \in F$ if

$$(33) \qquad \sup_{\theta \in \Theta, f \in F_{\delta_\varepsilon}(\bar{f})} \mathbf{E}_{\theta,f}[(\theta_\varepsilon^* - \theta)^2 I^\varepsilon(f)] = 1 + (1 + o(1)) \frac{r^\varepsilon}{\|\bar{f}'\|^2} \qquad \text{as } \varepsilon \to 0,$$

for some $\delta_\varepsilon > 0$ such that $\lim_{\varepsilon \to 0} \delta_\varepsilon = 0$.

Comparing Theorems 1 and 2 we see that if there exists a sequence of weights $h_k = \lambda_k^*$ for which Assumptions B and C are satisfied and

$$(34) \qquad \sup_{f \in F_{\delta_\varepsilon}(\bar{f})} R^\varepsilon[f, \lambda^*] \leq r^\varepsilon (1 + o(1)),$$

where $\lambda^* = \{\lambda_k^*\}$, then the estimator $\hat{\theta}_{\mathrm{AD}}$ with this choice of weights is second-order efficient. At first sight, it seems that one can take $\lambda_k^* = q_k$ from (25). However, for $h_k = q_k$ Assumption C is not fulfilled. Therefore we correct $q_k$, taking

$$\lambda_k^* = \begin{cases} 1, & k \leq \gamma_\varepsilon W_\varepsilon, \\ \left[1 - \left(\frac{k}{W_\varepsilon}\right)^{\beta-1}\right]_+, & k > \gamma_\varepsilon W_\varepsilon, \end{cases}$$

where $W_\varepsilon$ is a solution of (26) and $\gamma_\varepsilon = 1/\log(\varepsilon^{-2})$. For $k > \gamma_\varepsilon W_\varepsilon$, the weights $\lambda_k^*$ induce a prior on $\{f_k\}$ analogous to the one that appears in the proof of the lower bound of Theorem 2. The corresponding penalized maximum likelihood type estimator has the form

$$(35) \qquad \theta_{\mathrm{PML}}^* = \arg\max_{\tau \in \Theta} \left\{ \sum_{k=1}^\infty \lambda_k^* \left( \int_{-1/2}^{1/2} \cos[2\pi k(t - \tau)] x^\varepsilon(t) \, dt \right)^2 \right\}.$$

THEOREM 3. *Let a function $\bar{f} \in F$ be such that, for some $p > \beta > 1$,*

$$\sum_{k=1}^{\infty} (2\pi k)^{2p} \bar{f}_k^2 < \infty \tag{36}$$

*and $\lim_{\varepsilon \to 0} \delta_\varepsilon = 0$, $\lim_{\varepsilon \to 0} \delta_\varepsilon / (\varepsilon^2 W_\varepsilon^{1+\alpha}) = \infty$, for some $\alpha > 0$, where $W_\varepsilon$ satisfies (27). Then, as $\varepsilon \to 0$, the local asymptotic minimax risk admits the second-order expansion*

$$\inf_{\hat{\theta}_\varepsilon} \sup_{\theta \in \Theta,\, f \in F_{\delta_\varepsilon}(\bar{f})} \mathbf{E}_{\theta,f}[(\hat{\theta}_\varepsilon - \theta)^2 I^\varepsilon(f)] = 1 + (1+o(1)) \frac{r^\varepsilon}{\|\bar{f}'\|^2}. \tag{37}$$

*Moreover, the estimator $\theta^*_{\mathrm{PML}}$ defined in (35) is second-order efficient at $\bar{f}$.*

The proof of Theorem 3 is given in Section 9.

REMARK 1. Theorems 2 and 3 are local in $f$ and nonlocal in $\theta$. Inspection of the proofs shows that they can be turned into local ones in $\theta$ as well, that is, that one can replace $\sup_{\theta \in \Theta}$ by $\sup_{|\theta - \theta_0| \leq t}$ where $t > 0$ is a small number (fixed or tending to 0 with $\varepsilon$ not too fast) and $\theta_0$ is an interior point of $\Theta$.

REMARK 2. In the argument of Section 3, $\lambda_k = \sigma_k^2 / (\sigma_k^2 + \varepsilon^2)$. The values $(\sigma_k^*)^2$ corresponding to $\lambda_k^*$ for $k > \gamma_\varepsilon W_\varepsilon$ are thus

$$(\sigma_k^*)^2 = \frac{\lambda_k^* \varepsilon^2}{1 - \lambda_k^*} = \varepsilon^2 \left[ \left(\frac{W_\varepsilon}{k}\right)^{\beta-1} - 1 \right]_+.$$

One can interpret these $(\sigma_k^*)^2$ as variances of the prior distributions of the $f_k$'s introduced in Section 3. These variances appear also in the proof of the lower bound [cf. (46)]. The fact that the initial values of $\lambda_k^*$ are equal to 1 means that we do not put any prior distribution on the Fourier coefficients $f_k$ for $k \leq \gamma_\varepsilon W_\varepsilon$. Note that this is a particular choice of a prior associated with the Sobolev classes of functions.

REMARK 3. It is interesting to compare results on nonparametric and semiparametric second-order efficiency. Golubev and Levit [11, 12] and Dalalyan and Kutoyants [6] considered nonparametric problems where there exist $\sqrt{n}$-consistent first-order efficient estimators (such as estimation of the cumulative distribution function). In these problems there are simple efficient estimators, as the empirical c.d.f. and smoothed estimators allow one to improve upon these simple estimators, so that the second-order asymptotic terms are always *negative*. On the contrary, in semiparametric problems, as in the one considered here, simple empirical estimators are not efficient, and one has to use smoothing already to attain first-order efficiency. As we see from Theorems 1–3 (cf. also Golubev and Härdle [9, 10], who studied partial linear models), in semiparametric problems second-order asymptotic terms are *positive*, so that they always spoil asymptotics. This suggests that the choice of correct smoothing that allows one to optimize second-order asymptotic terms is more important in semiparametrics than in nonparametrics.

**7. Proof of Theorem 1.** In what follows we use the same notation $C$ for finite positive constants that may be different in different occasions and can depend only on $\tau_0, \rho, C_0, \rho_1$ and $C_1$.

The first step of the proof of Theorem 1 is to show that estimator $\hat{\theta}_{\mathrm{AD}}$ is $\varepsilon$-consistent.

7.1. *Consistency of $\hat{\theta}_{\mathrm{AD}}$.* The estimator $\hat{\theta}_{\mathrm{AD}}$ is a maximizer of the contrast function

$$L(\tau) = \sum_{k=1}^{\infty} h_k \left\{ \sqrt{2} \int_{-1/2}^{1/2} \cos[2\pi k(t-\tau)] x^{\varepsilon}(t) \, dt \right\}^2$$

$$= \sum_{k=1}^{\infty} h_k (f_k \cos[2\pi k(\tau - \theta)] + \varepsilon \xi_k(0) \cos(2\pi k \tau) + \varepsilon \xi_k^*(0) \sin(2\pi k \tau))^2$$

$$= \sum_{k=1}^{\infty} h_k f_k^2 \cos^2[2\pi k(\tau - \theta)] + 2\varepsilon \|f'\| \eta_1(\tau) + \varepsilon^2 \eta_2(\tau),$$

where $\theta$ is the true value of the parameter,

$$\xi_k(u) = \sqrt{2} \int_{-1/2}^{1/2} \cos[2\pi k(t-u)] n(t) \, dt,$$

$$\xi_k^*(u) = \sqrt{2} \int_{-1/2}^{1/2} \sin[2\pi k(t-u)] n(t) \, dt$$

and

$$\eta_1(\tau) = \frac{1}{\|f'\|} \sum_{k=1}^{\infty} h_k f_k \cos[2\pi k(\tau - \theta)] (\xi_k(0) \cos(2\pi k \tau) + \xi_k^*(0) \sin(2\pi k \tau)),$$

$$\eta_2(\tau) = \sum_{k=1}^{\infty} h_k (\xi_k(0) \cos(2\pi k \tau) + \xi_k^*(0) \sin(2\pi k \tau))^2.$$

The following three lemmas allow us to control the first derivatives of $\eta_1(\tau)$ and $\eta_2(\tau)$.

LEMMA 1. *Uniformly in $f \in F$ we have*

$$\mathbf{P}\left\{ \sup_{\tau \in \Theta} |\eta_1'(\tau)| > x \right\} \leq c_1 \exp(-c_2 x^2) \qquad \forall x > 0,$$

*where the constants $c_1 > 0$ and $c_2 > 0$ depend only on $\rho$ and $C_0$.*

PROOF. Note that $\eta_1'(\tau)$ is a stationary Gaussian random process with mean 0 and twice continuously differentiable correlation function $r(\cdot)$ such that $r''(0) \neq 0$. It follows from the Rice formula ([18], Theorem 7.3.2, or [5], Chapter 13.5, page 294; see also Proposition 2 in [28]) that for all $x > 0$,

$$(38) \quad \mathbf{P}\left\{\sup_{\tau \in \Theta} |\eta_1'(\tau)| > x\right\} \leq C\left[(r''(0)/r(0))^2 + 1\right] \exp\left(-\frac{x^2}{2r^2(0)}\right),$$

where $C > 0$ is a universal constant. Now, since $f \in F$,

$$r^2(0) = \mathbf{E}[\eta_1'(\tau)^2] = \|f'\|^{-2} \sum_{k=1}^{\infty} h_k^2 f_k^2 (2\pi k)^2 \geq (2\pi)^2 \|f'\|^{-2} \rho,$$

$$(r''(0))^2 = \mathbf{E}[\eta_1''(\tau)^2] = \|f'\|^{-2} \sum_{k=1}^{\infty} h_k^2 f_k^2 (2\pi k)^4 \leq C_0 \|f'\|^{-2},$$

which together with (38) proves the lemma. □

We will use the following simple fact about moderate deviations of the random variable:

$$\varsigma = \sum_{i=1}^{\infty} a_i (\xi_i^2 - 1),$$

where the $\xi_i$'s are i.i.d. standard normal random variables and $\{a_k\}$ is a sequence belonging to $\ell_2$, so that the random series converges almost surely.

LEMMA 2. *Let $a_k \not\equiv 0$, $\{a_k\} \in \ell_2$. For any $0 < x \leq \|a\|/\max_k |a_k|$ we have*

$$\mathbf{P}\{|\varsigma| \geq x\sqrt{\mathbf{E}[\varsigma^2]}\} \leq 2\exp(-x^2/16).$$

This result follows, for example, from (27) of Lemma 2 in [4].

LEMMA 3. *For any $0 < x \leq \|h'\|/\max_k h_k(2\pi k)$*

$$\mathbf{P}\left\{\sup_{\tau \in \Theta} |\eta_2'(\tau)| > 4\sum_{k=1}^{\infty} h_k(2\pi k) + x\|h'\|\right\} \leq 4\exp(-c_3 x^2),$$

*where $c_3 > 0$ is a universal constant.*

PROOF. Using the Cauchy–Schwarz inequality we get

$$\sup_\tau |\eta_2'(\tau)| \leq 2\sum_{k=1}^{\infty} h_k(2\pi k) \sup_\tau \{|\xi_k(0)\cos(2\pi k\tau) + \xi_k^*(0)\sin(2\pi k\tau)|$$
$$\times |-\xi_k(0)\sin(2\pi k\tau) + \xi_k^*(0)\cos(2\pi k\tau)|\}$$
$$\leq 2\sum_{k=1}^{\infty} h_k(2\pi k)[\xi_k^2(0) + \xi_k^{*2}(0)].$$

The rest follows from Lemma 2.  □

Consider now the expectation of the contrast function $L(\cdot)$

$$\mathbf{E}[L(\tau)] = \sum_{k=1}^{\infty} h_k f_k^2 \cos^2[2\pi k(\tau - \theta)].$$

LEMMA 4. *Let Assumptions* A1–A3 *and* B *be satisfied. Then*

$$\mathbf{E}[L(\tau)] - \mathbf{E}[L(\theta)] \leq -C|\tau - \theta|^2 \qquad \forall \tau \in \Theta,$$

*where the constant $C > 0$ depends only on $\tau_0$, $\rho$ and $C_0$.*

PROOF.  The derivatives of the function $G(\tau) = \mathbf{E}[L(\tau)]$ satisfy $G'(\theta) = 0$, $G''(\theta) = -2\sum_{k=1}^{\infty} h_k f_k^2 (2\pi k)^2 \leq -2(2\pi)^2 \rho$. Thus, the assertion of the lemma holds for $\tau$ in some neighborhood of $\theta$. Since also $\mathbf{E}[L(\tau)] < \mathbf{E}[L(\theta)]$ for all $\tau \in \Theta$, $\tau \neq \theta$, and $\Theta$ is a bounded interval (cf. Assumption A1), the lemma follows.  □

Now we are ready to show that $\hat{\theta}_{\text{AD}}$ is $\varepsilon$-consistent.

LEMMA 5.  *Let Assumptions* A1–A3 *and* B *be satisfied. Then, uniformly in $f \in F$ and in $\theta \in \Theta$,*

$$\mathbf{P}_{\theta,f}\{|\hat{\theta}_{\text{AD}} - \theta|\sqrt{I^\varepsilon(f)} > x\} \leq c_4 \exp(-c_5 x^2)$$

*for all $x \in [x_0, \|h'\|/\max_k h_k(2\pi k)]$, where $c_4 > 0, c_5 > 0, x_0 > 0$ are constants depending only on $\tau_0$, $\rho$, $C_0$, $C_1$.*

PROOF.  Due to Lemma 4 we have

$$\mathbf{P}_{\theta,f}\{|\hat{\theta}_{\text{AD}} - \theta|\sqrt{I^\varepsilon(f)} > x\}$$

$$\leq \mathbf{P}\Big\{\max_{\tau \in \Theta:\, |\tau-\theta|>x/\sqrt{I^\varepsilon(f)}} [L(\tau) - L(\theta)] \geq 0\Big\}$$

$$\leq \mathbf{P}\Big\{\max_{\tau \in \Theta:\, |\tau-\theta|>x/\sqrt{I^\varepsilon(f)}} [\mathbf{E}[L(\tau)] - \mathbf{E}[L(\theta)]$$

$$+ 2\varepsilon\|f'\|(\eta_1(\tau) - \eta_1(\theta))$$

$$+ \varepsilon^2(\eta_2(\tau) - \eta_2(\theta))] \geq 0\Big\}$$

$$\leq \mathbf{P}\Big\{\max_{\tau \in \Theta:\, |\tau-\theta|>x/\sqrt{I^\varepsilon(f)}} \Big[\mathbf{E}[L(\tau)] - \mathbf{E}[L(\theta)]$$

$$+ |\tau - \theta| \Big( 2\varepsilon \|f'\| \max_{t \in \Theta} |\eta_1'(t)|$$
$$+ \varepsilon^2 \max_{t \in \Theta} |\eta_2'(t)| \Big) \Big] \geq 0 \Big\}$$
$$\leq \mathbf{P} \Big\{ \max_{t \in \Theta} |\eta_1'(t)| + \varepsilon \max_{t \in \Theta} |\eta_2'(t)| \geq Cx \Big\}$$
$$\leq \mathbf{P} \Big\{ \max_{t \in \Theta} |\eta_1'(t)| \geq Cx \Big\} + \mathbf{P} \Big\{ \varepsilon \max_{t \in \Theta} |\eta_2'(t)| \geq Cx \Big\}.$$

The first probability on the last line is controlled by Lemma 1, whereas the second probability can be bounded, in view of Lemma 3, by $4\exp(-Cx^2)$, since according to Assumption B2 one has $\varepsilon \sum_{k=1}^{\infty} h_k(2\pi k) \leq C'$, $\varepsilon \|h'\| \leq C'$, where $C'$ depends only on $C_1$, and thus $Cx > 4\varepsilon \sum_{k=1}^{\infty} h_k(2\pi k) + c\varepsilon x \|h'\|$ for any $x \geq x_0$ if $x_0$ is large enough and $c > 0$ is small enough. □

7.2. *Proof of Theorem* 1. Let us introduce the event $\mathcal{A}_1 = \{|\hat{\theta}_{\mathrm{AD}} - \theta| \leq c_6 \varepsilon \sqrt{\log(\varepsilon^{-2})}\}$ where $c_6 > 0$ is a sufficiently large constant that can depend only on $\tau_0, \rho, C_0$ and $C_1$. The risk of $\hat{\theta}_{\mathrm{AD}}$ can be decomposed into two terms,

$$(39) \quad \mathbf{E}_{\theta,f}[(\hat{\theta}_{\mathrm{AD}} - \theta)^2] = \mathbf{E}_{\theta,f}[(\hat{\theta}_{\mathrm{AD}} - \theta)^2 \mathbb{1}_{\mathcal{A}_1}] + \mathbf{E}_{\theta,f}[(\hat{\theta}_{\mathrm{AD}} - \theta)^2 \mathbb{1}_{\mathcal{A}_1^c}].$$

Using (24) and Lemma 5 we find that, for $c_6$ large enough,

$$(40) \quad \mathbf{E}_{\theta,f}[(\hat{\theta}_{\mathrm{AD}} - \theta)^2 I^\varepsilon(f) \mathbb{1}_{\mathcal{A}_1^c}] \leq C\varepsilon^{-2} \mathbf{P}_{\theta,f}\{\mathcal{A}_1^c\} = O(\varepsilon^2) \qquad \text{as } \varepsilon \to 0.$$

Indeed, for $x = \varepsilon \sqrt{I^\varepsilon(f) \log(\varepsilon^{-2})} \geq C\sqrt{\log(\varepsilon^{-2})}$, due to Assumption B1 and (24), one has

$$\frac{x \max_k h_k(2\pi k)}{\|h'\|} \leq \frac{\sqrt{C_0} \log^{-3/2}(\varepsilon^{-2})}{\rho_1} \longrightarrow 0 \qquad \text{as } \varepsilon \to 0.$$

Thus we can apply Lemma 5, which yields (40) when $c_6$ is large enough. It remains to find the asymptotics of the first term on the right-hand side of (39). The estimator $\hat{\theta}_{\mathrm{AD}}$ satisfies

$$(41) \qquad L'(\hat{\theta}_{\mathrm{AD}}) = 0.$$

Using Taylor approximation of the left-hand side of (41) in a neighborhood of $\theta$ we may write, for some $\omega \in \Theta$,

$$(42) \qquad L_0(\theta) + (\theta - \hat{\theta}_{\mathrm{AD}}) L_1(\theta) + \tfrac{1}{2}(\theta - \hat{\theta}_{\mathrm{AD}})^2 L_2(\omega) = 0,$$

where

$$L_0(\theta) = \varepsilon \sum_{k=1}^{\infty} h_k(2\pi k) f_k \xi_k^*(\theta) + \varepsilon^2 \sum_{k=1}^{\infty} h_k(2\pi k) \xi_k^*(\theta) \xi_k(\theta),$$

$$L_1(\theta) = \sum_{k=1}^{\infty} h_k (2\pi k)^2 (f_k^2 + \varepsilon f_k \xi_k(\theta) + \varepsilon^2 [\xi_k^2(\theta) - \xi_k^{*2}(\theta)]),$$

$$L_2(\omega) = -8 \sum_{k=1}^{\infty} h_k (2\pi k)^3 \left( \int_{-1/2}^{1/2} \cos[2\pi k(t-\omega)] x^\varepsilon(t)\, dt \right)$$

$$\times \left( \int_{-1/2}^{1/2} \sin[2\pi k(t-\omega)] x^\varepsilon(t)\, dt \right).$$

LEMMA 6. *Let Assumptions* A1–A3 *and* B *be satisfied. Then*

$$\sup_{\theta \in \Theta, f \in F} \mathbf{E}_{\theta,f}[(L_1(\theta) - \mathbf{E}_{\theta,f}[L_1(\theta)])^2] = O(\varepsilon^2) \qquad \text{as } \varepsilon \to 0,$$

*and*

$$\sup_{\theta \in \Theta, f \in F} \mathbf{E}_{\theta,f}\left[ \sup_{\omega \in \Theta} |L_2(\omega)|^2 \right] \leq C.$$

PROOF. We omit the proof of the first relation since it follows from simple algebra. To prove the second one, using trigonometric formulae and the Cauchy–Schwarz inequality, we write

$$\sqrt{2} \left| \int \cos[2\pi k(t-\omega)] x^\varepsilon(t)\, dt \right|$$
$$= |f_k \cos[2\pi k(\theta - \omega)] + \varepsilon \xi_k(0) \cos[2\pi k\omega] + \varepsilon \xi_k^*(0) \sin[2\pi k\omega]|$$
$$\leq |f_k| + \varepsilon \sqrt{\xi_k(0)^2 + \xi_k^*(0)^2}.$$

Similarly, $\sqrt{2}|\int \sin[2\pi k(t-\omega)] x^\varepsilon(t)\, dt| \leq |f_k| + \varepsilon \sqrt{\xi_k(0)^2 + \xi_k^*(0)^2}$. Therefore

$$|L_2(\omega)| \leq C \sum_{k=1}^{\infty} h_k (2\pi k)^3 f_k^2 + C\varepsilon^2 \sum_{k=1}^{\infty} h_k (2\pi k)^3 [\xi_k^2(0) + \xi_k^{*2}(0)].$$

The second inequality of the lemma follows easily from this and Assumptions A3 and B2. □

To analyze the behavior of $\hat\theta_{\mathrm{AD}}$ we compare it to the root $\hat\tau$ of the linear equation

(43) $$L_0(\theta) + (\theta - \hat\tau) \mathbf{E}_{\theta,f}[L_1(\theta)] = 0$$

representing an approximation of (42).

LEMMA 7. *Let Assumptions* A1–A3, B *and* C *be satisfied. Then*

$$\mathbf{E}_{\theta,f}[(\hat\tau - \theta)^2 I^\varepsilon(f)] = 1 + (1 + o(1)) \frac{R^\varepsilon[f, h]}{\|f'\|^2},$$

*where* $o(1) \to 0$ *uniformly in* $f \in F$ *and in* $\theta \in \Theta$, *as* $\varepsilon \to 0$.

PROOF. Using the inequality $(1-h_k^2)f_k^2(2\pi k)^2 \leq 2(1-h_k)f_k^2(2\pi k)^2$, Assumption C and (24), we get from (43),

$$\mathbf{E}_{\theta,f}[(\hat{\tau} - \theta)^2 I^\varepsilon(f)]$$
$$= \frac{1 + \|f'\|^{-2}\sum_{k=1}^\infty [(h_k^2 - 1)f_k^2 + \varepsilon^2 h_k^2](2\pi k)^2}{[1 + \|f'\|^{-2}\sum_{k=1}^\infty (h_k - 1)f_k^2(2\pi k)^2]^2}$$
$$= \left[1 + \|f'\|^{-2}\sum_{k=1}^\infty [(h_k^2 - 1)f_k^2 + \varepsilon^2 h_k^2](2\pi k)^2\right]$$
$$\times \left[1 - 2\|f'\|^{-2}\sum_{k=1}^\infty (h_k - 1)f_k^2(2\pi k)^2 + o(R^\varepsilon[f,h])\right]$$
$$= 1 + (1+o(1))\|f'\|^{-2}R^\varepsilon[f,h]. \qquad \square$$

LEMMA 8. *Let Assumptions* A1–A3 *and* B *be satisfied. Then* $\mathbf{E}_{\theta,f}[(\hat{\theta}_{AD} - \hat{\tau})^2 \mathbb{1}_{\mathcal{A}_1}] \leq C\varepsilon^4 \log^2(\varepsilon^{-2})$.

PROOF. Since no confusion is possible, we omit the subscripts $\theta$, $f$ of the expectation. Subtracting (43) from (42) we obtain

$$(\hat{\theta}_{AD} - \hat{\tau})\mathbf{E}[L_1(\theta)] - (\theta - \hat{\theta}_{AD})(L_1(\theta) - \mathbf{E}[L_1(\theta)]) - \tfrac{1}{2}(\theta - \hat{\theta}_{AD})^2 L_2(\omega) = 0.$$

Note that $\mathbf{E}[L_1(\theta)] = \sum_k h_k(2\pi k)^2 f_k^2 \geq (2\pi)^2 \rho$ and that $(\hat{\theta}_{AD} - \theta)^2 \leq c_6^2 \varepsilon^2 \log(\varepsilon^{-2})$ on $\mathcal{A}_1$. Using these facts and Lemma 6 we get

$$\mathbf{E}[(\hat{\theta}_{AD} - \hat{\tau})^2 \mathbb{1}_{\mathcal{A}_1}]$$
$$\leq (\mathbf{E}[L_1(\theta)])^{-2}\bigg\{2\mathbf{E}[(\theta - \hat{\theta}_{AD})^2(L_1(\theta) - \mathbf{E}[L_1(\theta)])^2 \mathbb{1}_{\mathcal{A}_1}]$$
$$+ \mathbf{E}_{\theta,f}\bigg[(\theta - \hat{\theta}_{AD})^4 \sup_{\omega \in \Theta}|L_2(\omega)|^2 \mathbb{1}_{\mathcal{A}_1}\bigg]\bigg\}$$
$$\leq C\varepsilon^4 \log^2(\varepsilon^{-2}). \qquad \square$$

Now Assumption B1 and the fact that $h_1 = 1$ yield, for $\varepsilon$ small enough,

$$R^\varepsilon[f,h] \geq \varepsilon^2 \sum_{k=1}^\infty (2\pi k)^2 h_k^2 \geq \rho_1 \varepsilon^2 \left(\max_k h_k(2\pi k)\right)^2 \log^4(\varepsilon^{-2})$$
$$\geq \rho_1(2\pi)^2 \varepsilon^2 \log^4(\varepsilon^{-2}),$$

which implies that $\mathbf{E}_{\theta,f}[(\hat{\theta}_{AD} - \hat{\tau})^2 I^\varepsilon(f)\mathbb{1}_{\mathcal{A}_1}] = o(R^\varepsilon[f,h])$ uniformly in $f \in F$ and in $\theta \in \Theta$, as $\varepsilon \to 0$. This result together with (39), (40) and Lemma 7 completes the proof of Theorem 1.

**8. Proof of Theorem 2.** Before proceeding to the proof of Theorem 2 we give some preliminary results.

8.1. *An auxiliary Bayesian problem.* We consider a model with two observations that will be used as a building block for the subsequent proofs. Set

$$x = f_0 \cos(2\pi k\theta) + \varepsilon\xi, \qquad x^* = f_0 \sin(2\pi k\theta) + \varepsilon\xi^*,$$

where $\xi, \xi^*$ are independent $\mathcal{N}(0,1)$ random variables and $f_0$ is an $\mathcal{N}(\bar f, \sigma^2)$ random variable that does not depend on $(\xi, \xi^*)$, with $\bar f \in \mathbf{R}$, $\sigma^2 > 0$. Here $\theta$ is a parameter to be estimated based on the observations $x, x^*$ and $k$ is an integer. Define the Fisher information

$$\mathcal{J}_k^\varepsilon(\theta) = \mathbf{E}\left[\left(\frac{d}{d\theta} \log p_\theta(x,x^*)\right)^2\right],$$

where $p_\theta(x,x^*)$ is the probability density of the observations.

LEMMA 9. *We have* $\mathcal{J}_k^\varepsilon(\theta) = \varepsilon^{-2}(\bar f^2 + \frac{\sigma^4}{\varepsilon^2+\sigma^2})(2\pi k)^2$, *for any* $k \in \mathbb{Z}$.

PROOF. Denoting by $C$ multiplicative constants that do not depend on $\theta$, we have

$$p_\theta(x,x^*) = C \int \exp\Big\{-\frac{u^2}{2\sigma^2} - \frac{1}{2\varepsilon^2}[x - \bar f \cos(2\pi k\theta) - u\cos(2\pi k\theta)]^2$$

$$- \frac{1}{2\varepsilon^2}[x^* - \bar f \sin(2\pi k\theta) - u\sin(2\pi k\theta)]^2\Big\} du$$

$$= C \exp\Big\{\frac{\lambda}{2}[x\cos(2\pi k\theta) + x^*\sin(2\pi k\theta)]^2$$

$$+ (1-\lambda)\bar f[x\cos(2\pi k\theta) + x^*\sin(2\pi k\theta)]\Big\}$$

$$= C \exp\Big\{\frac{\lambda}{2\varepsilon^2}\Big[x\cos(2\pi k\theta) + x^*\sin(2\pi k\theta) + \frac{1-\lambda}{\lambda}\bar f\Big]^2\Big\},$$

where $\lambda = \sigma^2/(\varepsilon^2+\sigma^2)$. Hence writing $f_0 = \bar f + \eta\sigma$ where $\eta \sim \mathcal{N}(0,1)$ and $\eta$ is independent of $(\xi,\xi^*)$, one obtains

$$\mathcal{J}_k^\varepsilon(\theta) = \mathbf{E}\left[\left(\frac{\lambda}{2\varepsilon^2}\frac{d}{d\theta}\Big[x\cos(2\pi k\theta) + x^*\sin(2\pi k\theta) + \frac{1-\lambda}{\lambda}\bar f\Big]^2\right)^2\right]$$

$$= (2\pi k)^2 \varepsilon^{-4} \lambda^2 \mathbf{E}\bigg[\Big(\bar f + \eta\sigma + \varepsilon\xi\cos(2\pi k\theta) + \varepsilon\xi^*\sin(2\pi k\theta) + \frac{1-\lambda}{\lambda}\bar f\Big)^2$$

$$\times (-\varepsilon\xi\sin(2\pi k\theta) + \varepsilon\xi^*\cos(2\pi k\theta))^2\bigg]$$

$$= (2\pi k)^2 \varepsilon^{-2} \lambda^2 \mathbf{E}\bigg[\bigg((\lambda^{-1}\bar{f} + \eta\sigma)(-\xi\sin(2\pi k\theta) + \xi^*\cos(2\pi k\theta))$$
$$+ \frac{\varepsilon}{2}(\xi^{*2} - \xi^2)\sin(4\pi k\theta) + \varepsilon\xi\xi^*\cos(4\pi k\theta)\bigg)^2\bigg]$$
$$= (2\pi k)^2 \varepsilon^{-2} \lambda^2 [\lambda^{-2}\bar{f}^2 + \sigma^2 + \varepsilon^2]. \qquad \square$$

8.2. *Lower bounds for Bayes risks.* In this subsection we consider the sequence model (17) where we suppose that the $f_k$'s are no longer fixed values but independent random variables distributed as $\mathcal{N}(\bar{f}_k, \sigma_k^2)$ with some $\sigma_k \geq 0$. By convention, $\sigma_k = 0$ means that the corresponding $f_k$ is equal to $\bar{f}_k$ almost surely. We assume in what follows that $\sigma_k > 0$ only for a finite (and possibly depending on $\varepsilon$) number of indices $k$. We also assume that the random sequence $(f_k, k = 1, 2, \ldots)$ does not depend on the noises $(\xi_k, \xi_k^*, k = 1, 2, \ldots)$. We will refer to this model as *the Bayes model with fixed $\theta$*. Let $\Psi_\sigma(df)$ denote the probability distribution of $f = \{f_k\} \in \ell_2$ in this model.

Along with this, we will consider *the full Bayes model* defined in the same way, except that in this new model $\theta$ is supposed to be a random variable having a density $\pi(x)$, $x \in \Theta$, that vanishes at the endpoints of the interval $\Theta$ and has finite Fisher information $I_\pi = \int (\pi'(x))^2 \pi^{-1}(x)\,dx$. It will be assumed that $\theta$ is independent of $(f_k, \xi_k, \xi_k^*, k = 1, 2, \ldots)$.

We denote by $\mathbb{E}$ the expectation with respect to the joint distribution of $(x_k, x_k^*, k = 1, 2, \ldots)$ and $\theta$ in the full Bayes model and by $\mathbb{E}_\theta$ the expectation w.r.t. the distribution of $(x_k, x_k^*, k = 1, 2, \ldots)$ in the Bayes model with fixed $\theta$. Define
$$\lambda_k = \frac{\sigma_k^2}{\varepsilon^2 + \sigma_k^2}, \qquad k = 1, 2, \ldots.$$

LEMMA 10. *Assume that the density $\pi(x)$ vanishes at the endpoints of the interval $\Theta$ and has finite Fisher information $I_\pi$. Then*

$$(44) \qquad \inf_{\hat{\theta}_\varepsilon} \mathbb{E}[(\hat{\theta}_\varepsilon - \theta)^2 \bar{I}^\varepsilon] \geq 1 + \frac{1}{\bar{I}^\varepsilon}\sum_{k=1}^\infty (2\pi k)^2 \lambda_k + O(\varepsilon^2),$$

*where*
$$\bar{I}^\varepsilon = \int I^\varepsilon(f)\, \Psi_\sigma(df) = \varepsilon^{-2} \sum_{k=1}^\infty (2\pi k)^2 (\bar{f}_k^2 + \sigma_k^2).$$

The proofs of this and subsequent lemmas are given in the Appendix.

In the next subsection we will show that one can choose the sequence $\{\sigma_k\}$ so that the right-hand side of (44) coincides asymptotically with the lower bound of Theorem 2. However, the left-hand side of (44) is different from that of (32). One difference is that in (32) the risk is normalized by the Fisher information $I^\varepsilon(f)$, while in (44) we have its average $\bar I^\varepsilon$ w.r.t. the distribution $\Psi_\sigma(df)$. The next lemma shows that $I^\varepsilon(f)$ is sufficiently close to $\bar I^\varepsilon$; in particular, its variance is small enough.

LEMMA 11. *If the $\sigma_k^2$'s are such that*

$$\text{(45)} \qquad \sum_{k=1}^\infty (2\pi k)^4 \sigma_k^4 + \sup_k \sigma_k^2 = o\left(\varepsilon^2 \sum_{k=1}^\infty (2\pi k)^2 \lambda_k\right),$$

*then*

$$\int (I^\varepsilon(f) - \bar I^\varepsilon)^2 \Psi_\sigma(df) = o\left(\varepsilon^{-2} \sum_{k=1}^\infty (2\pi k)^2 \lambda_k\right) \qquad as\ \varepsilon \to 0.$$

LEMMA 12. *Assume that the density $\pi(x)$ vanishes at the endpoints of the interval $\Theta$ and has finite Fisher information $I_\pi = \int (\pi'(x))^2 \pi^{-1}(x)\,dx$. Then, for any $f \in F$,*

$$\inf_{\hat\theta_\varepsilon} \int_\Theta \mathbf{E}_{\theta,f}[(\hat\theta_\varepsilon - \theta)^2 I^\varepsilon(f)]\, \pi(\theta)\,d\theta \geq \frac{I^\varepsilon(f)}{I^\varepsilon(f) + I_\pi} \geq 1 - \frac{I_\pi}{I^\varepsilon(f)}.$$

Proof of this lemma is omitted: this is the standard Van Trees inequality for the problem of estimation of $\theta$ with fixed $f$ in model (1) ([33]; see also [7]).

LEMMA 13. *If the sequence $\{\sigma_k\}$ satisfies relation (45) and $\delta < \sqrt{\rho/2}$, then*

$$\int_{F_\delta(\bar f)} \left(1 - \frac{\bar I^\varepsilon}{I^\varepsilon(f)}\right) \Psi_\sigma(df) \leq o\left(\varepsilon^2 \sum_{k=1}^\infty (2\pi k)^2 \lambda_k\right) + C\mathbf{P}(f \notin F_\delta(\bar f)).$$

8.3. *From Bayes to minimax bounds.* The main idea of the proof of Theorem 2 is to bound from below the minimax risk by a suitably chosen Bayes risk. In the rest of this section we consider the full Bayes model defined in Section 8.2 with a special choice of the $\sigma_k$'s. Namely, we set

$$\text{(46)} \qquad \sigma_k^2 = \begin{cases} 0, & k \leq \gamma_\varepsilon W_\varepsilon, \\ (1 - \gamma_\varepsilon) s_k^2, & k > \gamma_\varepsilon W_\varepsilon, \end{cases}$$

where $W_\varepsilon$ is a solution of (26), $s_k^2$ is defined by (28) and $\gamma_\varepsilon = 1/\log(\varepsilon^{-2})$ (here and later we suppose that $\varepsilon$ is small enough, so that $\gamma_\varepsilon < 1$). To derive the minimax lower bound of Theorem 2 from the Bayes bounds of Section 8.2 we need first to show that with a probability close to 1 the Gaussian random sequence $\{f_k\}$ belongs to the set $F_\delta(\bar f)$. In fact, the following result holds.

LEMMA 14. *For any $\delta^2 \geq \varepsilon^2 W_\varepsilon \gamma_\varepsilon^{2-2\beta}$ we have $\mathbf{P}\{f \notin F_\delta(\bar{f})\} \leq e^{-C\gamma_\varepsilon^2 W_\varepsilon}$.*

8.4. *Proof of Theorem* 2. Recall that we consider the full Bayes model with the $\sigma_k^2$'s chosen according to (46) and $\lambda_k = \sigma_k^2/(\varepsilon^2 + \sigma_k^2)$. Note that in this case

$$\text{(47)} \qquad \varepsilon^2 \sum_{k=1}^{\infty} (2\pi k)^2 \lambda_k = r^\varepsilon(1 + o(1)) \qquad \text{as } \varepsilon \to 0.$$

Indeed, (25) and (28) imply that $|\lambda_k/q_k - 1| \leq \gamma_\varepsilon$ for $k > \gamma_\varepsilon W_\varepsilon$, and hence [cf. (29)]

$$\varepsilon^2 \sum_{k=1}^{\infty} (2\pi k)^2 \lambda_k = (1 + o(1))\varepsilon^2 \sum_{k > \gamma_\varepsilon W_\varepsilon} (2\pi k)^2 q_k$$

$$= (1 + o(1))\left(r^\varepsilon - \varepsilon^2 \sum_{k \leq \gamma_\varepsilon W_\varepsilon} (2\pi k)^2 q_k\right).$$

Here [cf. (30)]

$$\varepsilon^2 \sum_{k \leq \gamma_\varepsilon W_\varepsilon} (2\pi k)^2 q_k = (2\pi)^2 \varepsilon^2 W_\varepsilon^3 \sum_{k \leq \gamma_\varepsilon W_\varepsilon} \left(\frac{k}{W_\varepsilon}\right)^2 \left[1 - \left(\frac{k}{W_\varepsilon}\right)^{\beta-1}\right] \frac{1}{W_\varepsilon}$$

$$\leq C\varepsilon^2 W_\varepsilon^3 \int_0^{\gamma_\varepsilon} (x^2 + x^{\beta+1})\, dx \leq C\gamma_\varepsilon^3 \varepsilon^2 W_\varepsilon^3 = o(r^\varepsilon),$$

and thus (47) follows.

Next, we check that if the $\sigma_k^2$'s are chosen according to (46), then condition (45) is satisfied, so that one can apply Lemmas 10–13. In fact, (27) yields $W_\varepsilon \asymp \varepsilon^{-2/(2\beta+1)}$ with $\beta > 1$, and using (47), (28) and (30) we get, as $\varepsilon \to 0$,

$$\varepsilon^2 \sum_{k=1}^{\infty} (2\pi k)^2 \lambda_k \asymp \varepsilon^2 W_\varepsilon^3 \to 0,$$

$$\sum_{k=1}^{\infty} (2\pi k)^4 \sigma_k^4 \leq \varepsilon^4 \sum_{\gamma_\varepsilon W_\varepsilon \leq k \leq W_\varepsilon} (2\pi k)^4 (W_\varepsilon/k)^{2(\beta-1)} \leq C\varepsilon^4 W_\varepsilon^5 \gamma_\varepsilon^{2-2\beta}$$

$$= o\left(\varepsilon^2 \sum_{k=1}^{\infty} (2\pi k)^2 \lambda_k\right),$$

$$\sup_k \sigma_k^2 \leq \varepsilon^2 \gamma_\varepsilon^{1-\beta} = o\left(\varepsilon^2 \sum_{k=1}^{\infty} (2\pi k)^2 \lambda_k\right).$$

Now we start the main body of the proof of Theorem 2. First note that, in a standard way, conditioning on $(x_k, x_k^*,\ k = 1, 2, \dots)$ and using Jensen's inequality, one can easily show that it is sufficient to prove the lower bound of Theorem 2 for estimators $\hat\theta_\varepsilon$ depending on $\mathbf{X}^\varepsilon$ only via $(x_k, x_k^*,\ k=1,2,\dots)$. Let $\mathcal{T}_\varepsilon$ denote the set of all estimators $\hat\theta_\varepsilon$ of $\theta$ measurable with respect to $(x_k, x_k^*,\ k=1,2,\dots)$ and satisfying the inequalities

$$(48)\qquad \sup_{f\in F_{\delta_\varepsilon}(\bar f)}\sup_{\theta\in\Theta}\mathbf{E}_{\theta,f}[(\hat\theta_\varepsilon-\theta)^2 I^\varepsilon(f)]\leq 1+\frac{2r^\varepsilon}{\|\bar f'\|^2}\quad\text{and}\quad |\hat\theta_\varepsilon|\leq 1.$$

It is enough to restrict our attention to the estimators from $\mathcal{T}_\varepsilon$, since for estimators that do not satisfy one of the inequalities in (48) the lower bound of Theorem 2 is evident.

Clearly,

$$(49)\qquad \sup_{f\in F_{\delta_\varepsilon}(\bar f)}\int_\Theta \mathbf{E}_{\theta,f}[(\hat\theta_\varepsilon-\theta)^2 I^\varepsilon(f)]\pi(\theta)\,d\theta \leq 1+\frac{2r^\varepsilon}{\|\bar f'\|^2}\qquad \forall\hat\theta_\varepsilon\in\mathcal{T}_\varepsilon.$$

We have

$$\inf_{\hat\theta\in\mathcal{T}_\varepsilon}\sup_{\theta,\,f\in F_\delta(\bar f)}\mathbf{E}_{\theta,f}[(\hat\theta-\theta)^2 I^\varepsilon(f)]$$

$$\geq \inf_{\hat\theta\in\mathcal{T}_\varepsilon}\mathbb{E}[(\hat\theta-\theta)^2 I^\varepsilon(f)\mathbb{1}_{F_\delta(\bar f)}(f)]$$

$$(50)\qquad \geq \inf_{\hat\theta\in\mathcal{T}_\varepsilon}\mathbb{E}[(\hat\theta-\theta)^2 \bar I^\varepsilon \mathbb{1}_{F_\delta(\bar f)}(f)] - \sup_{\hat\theta\in\mathcal{T}_\varepsilon}\mathbb{E}[(\hat\theta-\theta)^2(\bar I^\varepsilon - I^\varepsilon(f))\mathbb{1}_{F_\delta(\bar f)}]$$

$$\geq \inf_{\hat\theta\in\mathcal{T}_\varepsilon}\mathbb{E}[(\hat\theta-\theta)^2 \bar I^\varepsilon] - o(\varepsilon^2) - \sup_{\hat\theta\in\mathcal{T}_\varepsilon}\mathbb{E}[(\hat\theta-\theta)^2(\bar I^\varepsilon - I^\varepsilon(f))\mathbb{1}_{F_\delta(\bar f)}(f)]$$

$$\geq \inf_{\hat\theta}\mathbb{E}[(\hat\theta-\theta)^2 \bar I^\varepsilon] - o(\varepsilon^2) - \sup_{\hat\theta\in\mathcal{T}_\varepsilon}\mathbb{E}[(\hat\theta-\theta)^2(\bar I^\varepsilon - I^\varepsilon(f))\mathbb{1}_{F_\delta(\bar f)}(f)],$$

where we have used the inequality

$$\sup_{\hat\theta\in\mathcal{T}_\varepsilon}\mathbb{E}[(\hat\theta-\theta)^2 \bar I^\varepsilon \mathbb{1}_{F_\delta^c(\bar f)}(f)] \leq C\varepsilon^{-2}\exp(-C\gamma^2 W_\varepsilon) = o(\varepsilon^2),$$

which is a direct consequence of the estimates $|\hat\theta|\leq 1$, $|\theta|<1/4$, $\bar I^\varepsilon \leq C\varepsilon^{-2}$, relation (27) and Lemma 14. The last term in (50) can be represented as

$$\mathbb{E}[(\hat\theta-\theta)^2(\bar I^\varepsilon - I^\varepsilon(f))\mathbb{1}_{F_\delta(\bar f)}(f)]$$

$$(51)\qquad = \int_{F_\delta(\bar f)}\left(1-\frac{\bar I^\varepsilon}{I^\varepsilon(f)}\right)\Psi_\sigma(df)$$

$$\qquad + \mathbb{E}([(\hat\theta-\theta)^2 I^\varepsilon(f)-1](1-\bar I^\varepsilon/I^\varepsilon(f))\mathbb{1}_{F_\delta(\bar f)}(f)).$$

Due to Lemmas 13 and 14, the second term on the right-hand side of (51) is asymptotically negligible with respect to $\varepsilon^2\sum_k(2\pi k)^2\lambda_k = r^\varepsilon(1+o(1))$ [cf.

(47)]. To evaluate the first term, note that

$$\mathbb{E}([(\hat{\theta}-\theta)^2 I^\varepsilon(f) - 1](1 - \bar{I}^\varepsilon/I^\varepsilon(f))\mathbb{1}_{F_\delta(\bar{f})}(f))$$

$$(52) \quad \leq \sup_{f \in F_\delta(\bar{f})} \left| \int_\Theta \mathbf{E}_{\theta,f}[(\hat{\theta}-\theta)^2 I^\varepsilon(f) - 1]\pi(\theta)\,d\theta \right| \mathbb{E}(|1 - \bar{I}^\varepsilon/I^\varepsilon(f)|)$$

$$\leq C\varepsilon^2 \sup_{f \in F_\delta(\bar{f})} \left| \int_\Theta \mathbf{E}_{\theta,f}[(\hat{\theta}-\theta)^2 I^\varepsilon(f)]\pi(\theta)\,d\theta - 1 \right| [\mathbb{E}(\bar{I}^\varepsilon - I^\varepsilon(f))^2]^{1/2}.$$

It follows from (58) and Lemma 11 that $\varepsilon^4 \mathbb{E}(\bar{I}^\varepsilon - I^\varepsilon(f))^2$ is $o(1)$. Now, Lemma 12, inequality (49) and the fact that $\sup_{f \in F_\delta(\bar{f})} I_\pi/I^\varepsilon(f) \leq C\varepsilon^2 = o(r^\varepsilon)$ [cf. (58)] imply that

$$(53) \quad \sup_{f \in F_\delta(\bar{f})} \left| \int_\Theta \mathbf{E}_{\theta,f}[(\hat{\theta}_\varepsilon - \theta)^2 I^\varepsilon(f)]\pi(\theta)\,d\theta - 1 \right| \leq Cr^\varepsilon \qquad \forall \hat{\theta}_\varepsilon \in \mathcal{T}_\varepsilon.$$

Plugging (51)–(53) in (50) and using Lemmas 10, 13 and 14, we get

$$\inf_{\hat{\theta} \in \mathcal{T}_\varepsilon} \sup_{\theta \in \Theta, f \in F_\delta(\bar{f})} \mathbf{E}_{\theta,f}[(\hat{\theta}-\theta)^2 I^\varepsilon(f)] \geq \inf_{\hat{\theta}} \mathbb{E}[(\hat{\theta}-\theta)^2 \bar{I}^\varepsilon] + o(r^\varepsilon)$$

$$\geq 1 + \frac{1}{\bar{I}^\varepsilon} \sum_{k=1}^\infty (2\pi k)^2 \lambda_k + o(r^\varepsilon)$$

$$= 1 + \frac{r^\varepsilon}{\|\bar{f}'\|^2} + o(r^\varepsilon),$$

where for the last equality we have used (47) and the fact that, due to (28) and (46),

$$|\varepsilon^2 \bar{I}^\varepsilon - \|\bar{f}'\|^2| \leq \sum_{k > \gamma_\varepsilon W_\varepsilon} (2\pi k)^2 \sigma_k^2$$

$$< \varepsilon^2 \sum_{\gamma_\varepsilon W_\varepsilon \leq k \leq W_\varepsilon} (2\pi k)^2 (W_\varepsilon/k)^{\beta-1}$$

$$\leq C\varepsilon^2 W_\varepsilon^3 \gamma_\varepsilon^{1-\beta} = o(1),$$

as $\varepsilon \to 0$.

**9. Proof of Theorem 3.** It is enough to check that Assumptions B and C are satisfied for $h_k = \lambda_k^*$ and that (34) holds. We first check Assumption C. Recall that we supposed w.l.o.g. that $\gamma_\varepsilon < 1$. Then $1 - \lambda_k^* \geq \gamma_\varepsilon^{\beta-1}$ for $k > \gamma_\varepsilon W_\varepsilon$, and we have

$$\sum_{k=1}^\infty (1-\lambda_k^*)(2\pi k)^2 f_k^2 = \sum_{k > \gamma_\varepsilon W_\varepsilon} (1-\lambda_k^*)(2\pi k)^2 f_k^2$$

$$\leq \gamma_\varepsilon^{1-\beta} \sum_{k=1}^\infty (1-\lambda_k^*)^2 (2\pi k)^2 f_k^2 \leq \gamma_\varepsilon^{1-\beta} r^\varepsilon.$$

This and (30) show that Assumption C is satisfied for $h_k = \lambda_k^*$. Using (27) we find that Assumption B also holds. Indeed, Assumption B2 amounts to checking that $\varepsilon^2 W_\varepsilon^5 \leq C_1$, which is clearly the case for $\beta > 1$, whereas Assumption B1 follows from the relation $\sqrt{W_\varepsilon}/\log^2(\varepsilon^{-2}) \to +\infty$, as $\varepsilon \to 0$. Now we are ready to check (34). For any $\kappa \in [0,1]$ one obtains [recall that the $q_k$'s are defined by (28)]

$$
\begin{aligned}
\sup_{f\in F_{\delta_\varepsilon}(\bar f)} R^\varepsilon[f,\lambda^*] &\leq \sup_{v\in \mathcal{W}(\beta,L)} \sum_{k=1}^\infty (1-\lambda_k^*)^2 (2\pi k)^2 (\bar f_k + v_k)^2 \\
&\quad + \sum_{k=1}^\infty (2\varepsilon\pi k \lambda_k^*)^2 \\
&\leq (1+\kappa) r^\varepsilon + \frac{2}{\kappa} \sum_{k=1}^\infty (1-\lambda_k^*)^2 (2\pi k)^2 \bar f_k^2 \\
&\quad + \varepsilon^2 \sum_{k\leq \gamma_\varepsilon W_\varepsilon} (1-q_k^2)(2\pi k)^2,
\end{aligned}
$$
(54)

where for the last inequality we have used (29). From (36) we obtain

$$\sum_{k=1}^\infty (1-\lambda_k^*)^2 (2\pi k)^2 \bar f_k^2 \leq \sum_{k>\gamma_\varepsilon W_\varepsilon} (2\pi k)^2 \bar f_k^2 \leq C(\gamma_\varepsilon W_\varepsilon)^{2-2p} = o(r^\varepsilon).$$

Note also that due to the relations $r^\varepsilon \geq C\varepsilon^2 W_\varepsilon^3$ [cf. (30)] and $\gamma_\varepsilon \to 0$ we get

$$\varepsilon^2 \sum_{k\leq \gamma_\varepsilon W_\varepsilon} (q_k^2-1)(2\pi k)^2 \leq 2\varepsilon^2 \sum_{k\leq \gamma_\varepsilon W_\varepsilon} \left(\frac{k}{W_\varepsilon}\right)^{\beta-1} (2\pi k)^2 \leq C\varepsilon^2 \gamma_\varepsilon^{\beta+2} W_\varepsilon^3 = o(r^\varepsilon).$$

These inequalities and (54) prove (34), since $\kappa$ can be arbitrarily small.

## APPENDIX

PROOF OF LEMMA 10. We start by applying the Van Trees inequality ([33]; see also [7]):

(55) $$\inf_{\hat\theta_\varepsilon} \mathbb{E}[(\hat\theta_\varepsilon - \theta)^2] \geq \left(\int_\Theta \mathcal{J}^\varepsilon(\theta)\pi(\theta)\,d\theta + I_\pi\right)^{-1},$$

where $\mathcal{J}^\varepsilon(\theta)$ is the Fisher information on $\theta$ contained in the observations $(x_k, x_k^*, k=1,2,\dots)$ for the Bayes model with fixed $\theta$. Since these observations are independent, $\mathcal{J}^\varepsilon(\theta)$ is the sum over $k$ of the Fisher information of pairs $(x_k, x_k^*)$. So using Lemma 9 we get that $\mathcal{J}^\varepsilon(\theta)$ does not depend on $\theta$ and equals

$$\mathcal{J}^\varepsilon(\theta) = \varepsilon^{-2}\sum_{k=1}^{\infty}(2\pi k)^2(\bar{f}_k^2 + \lambda_k \sigma_k^2) = \bar{I}^\varepsilon - \sum_{k=1}^{\infty}(2\pi k)^2 \lambda_k.$$

Therefore

(56) $$\bar{I}^\varepsilon \left( \int_\Theta \mathcal{J}^\varepsilon(\theta)\pi(\theta)\,d\theta + I_\pi \right)^{-1} \geq 1 + \frac{1}{\bar{I}^\varepsilon}\sum_{k=1}^{\infty}(2\pi k)^2 \lambda_k - \frac{I_\pi}{\bar{I}^\varepsilon}.$$

To complete the proof, it is enough to remark that, in view of (24),

(57) $$\bar{I}^\varepsilon \geq \varepsilon^{-2}\|\bar{f}'\|^2 \geq \varepsilon^{-2}(2\pi)^2 \rho. \qquad \square$$

PROOF OF LEMMA 11. Using the independence of the $f_k$'s for different values of $k$, we get

$$\varepsilon^4 \int (I^\varepsilon(f) - \bar{I}^\varepsilon)^2 \Psi_\sigma(df) = \sum_{k=1}^{\infty}(2\pi k)^4 \int (f_k^2 - \bar{f}_k^2 - \sigma_k^2)^2 \Psi_\sigma(df)$$

$$= \sum_{k=1}^{\infty}(2\pi k)^4(4\bar{f}_k^2\sigma_k^2 + 2\sigma_k^4)$$

$$\leq 4\left[\|\bar{f}''\|^2 \sup_k \sigma_k^2 + \sum_{k=1}^{\infty}(2\pi k)^4 \sigma_k^4\right].$$

The assertion of the lemma follows now from (45). $\square$

PROOF OF LEMMA 13. First note that using Assumption A2 one obtains

(58) $$\varepsilon^2 I^\varepsilon(f) \geq (2\pi)^2 f_1^2 \geq (2\pi)^2(\bar{f}_1^2 - \delta^2) > (2\pi)^2 \rho/2,$$

for any $f \in F_\delta(\bar{f})$ and for any $\delta < \sqrt{\rho/2}$. Furthermore, by (24),

(59) $$\varepsilon^2 I^\varepsilon(f) \leq 2(\|\bar{f}'\|^2 + L) \leq 2(C_0 + L) \qquad \forall\, f \in F_\delta(\bar{f}).$$

The elementary identity $1 - y = y^{-1} - 1 - y(1 - y^{-1})^2$ yields

$$\int_{F_\delta(\bar{f})}\left(1 - \frac{\bar{I}^\varepsilon}{I^\varepsilon(f)}\right)\Psi_\sigma(df) = \int_{F_\delta(\bar{f})}\left(\frac{I^\varepsilon(f)}{\bar{I}^\varepsilon} - 1\right)\Psi_\sigma(df)$$

$$+ \int_{F_\delta(\bar{f})}\left(\frac{I^\varepsilon(f)}{\bar{I}^\varepsilon} - 1\right)^2 \frac{\bar{I}^\varepsilon}{I^\varepsilon(f)}\Psi_\sigma(df).$$

To estimate the first integral on the right-hand side we note that $\bar{I}^\varepsilon = \int I^\varepsilon(f) \times \Psi_\sigma(df)$; therefore using (57) and (59) we get

$$\left|\int_{F_\delta(\bar{f})}\left(\frac{I^\varepsilon(f)}{\bar{I}^\varepsilon} - 1\right)\Psi_\sigma(df)\right| = \left|\int_{F_\delta^c(\bar{f})}\left(\frac{I^\varepsilon(f)}{\bar{I}^\varepsilon} - 1\right)\Psi_\sigma(df)\right|$$

$$\leq C\mathbf{P}(f \notin F_\delta(\bar{f})),$$

where $F_\delta^c(\bar{f}) = \ell^2 \setminus F_\delta(\bar{f})$. Finally, due to (57) and (58),

$$\int_{F_\delta(\bar{f})} \left(\frac{I^\varepsilon(f)}{\bar{I}^\varepsilon} - 1\right)^2 \frac{\bar{I}^\varepsilon}{I^\varepsilon(f)} \Psi_\sigma(df) \leq C\varepsilon^4 \int (I^\varepsilon(f) - \bar{I}^\varepsilon)^2 \Psi_\sigma(df).$$

The rest follows from Lemma 11. □

PROOF OF LEMMA 14. Let $\eta_k$ be i.i.d. $\mathcal{N}(0,1)$ random variables. We have

(60)
$$\mathbf{P}\{f \notin F_\delta(\bar{f})\} \leq \mathbf{P}\left\{\sum_{k > \gamma_\varepsilon W_\varepsilon} \eta_k^2 s_k^2 > \delta^2\right\}$$
$$+ \mathbf{P}\left\{\sum_{k > \gamma_\varepsilon W_\varepsilon} (2\pi k)^{2\beta} \eta_k^2 s_k^2 > \frac{L}{1 - \gamma_\varepsilon}\right\}.$$

We use Lemma 2 in order to evaluate the second probability. Note that

$$\sum_{k > \gamma_\varepsilon W_\varepsilon} (2\pi k)^{4\beta} s_k^4 \leq C\varepsilon^4 W_\varepsilon^{4\beta+1}$$

and $\max_k s_k^2 (2\pi k)^{2\beta} \leq C\varepsilon^2 W_\varepsilon^{2\beta}$. Therefore, by Lemma 2, for any $x \leq C\sqrt{W_\varepsilon}$ we have

$$\mathbf{P}\left\{\sum_{k \geq \gamma_\varepsilon W_\varepsilon} (2\pi k)^{2\beta} (\eta_k^2 - 1) s_k^2 > x\varepsilon^2 W_\varepsilon^{2\beta+1/2}\right\} \leq \exp(-Cx^2).$$

Applying this inequality for $x = \gamma_\varepsilon L / \varepsilon^2 W_\varepsilon^{2\beta+1/2}$ [note that in view of (27) $x$ is less than $C\sqrt{W_\varepsilon}$], and using the fact that $\sum_{k > \gamma_\varepsilon W_\varepsilon} (2\pi k)^{2\beta} s_k^2 \leq C\varepsilon^2 W_\varepsilon^{2\beta} \gamma_\varepsilon^{1-\beta} = o(\gamma_\varepsilon)$, one obtains

$$\mathbf{P}\left\{\sum_{k > \gamma_\varepsilon W_\varepsilon} (2\pi k)^{2\beta} \eta_k^2 s_k^2 > \frac{L}{(1 - \gamma_\varepsilon)}\right\} \leq \mathbf{P}\left\{\sum_{k > \gamma_\varepsilon W_\varepsilon} (2\pi k)^{2\beta}(\eta_k^2 - 1) s_k^2 > \frac{\gamma_\varepsilon L}{2}\right\}$$
$$\leq \exp\left\{-\frac{C\gamma_\varepsilon^2 L^2}{\varepsilon^4 W_\varepsilon^{4\beta+1}}\right\} \leq \exp(-C\gamma_\varepsilon^2 W_\varepsilon).$$

The first probability on the right-hand side of (60) can be estimated similarly. We have $\sum_{k > \gamma_\varepsilon W_\varepsilon} s_k^4 \leq C\varepsilon^4 W_\varepsilon \gamma_\varepsilon^{-2\beta+3}$ and $\max_{k > \gamma W_\varepsilon} s_k^2 \leq C\varepsilon^2 \gamma_\varepsilon^{-\beta+1}$. Hence, by Lemma 2, for any $x \leq C\sqrt{\gamma_\varepsilon W_\varepsilon}$,

$$\mathbf{P}\left\{\sum_{k > \gamma_\varepsilon W_\varepsilon} (\eta_k^2 - 1) s_k^2 > x\varepsilon^2 \sqrt{W_\varepsilon} \gamma_\varepsilon^{-2\beta+3/2}\right\} \leq \exp(-Cx^2).$$

So with $x = C\sqrt{\gamma_\varepsilon W_\varepsilon}$, noting that $\sum_{k > \gamma_\varepsilon W_\varepsilon} s_k^2 \leq \varepsilon^2 W_\varepsilon \gamma_\varepsilon^{2-\beta}$, one obtains

$$\mathbf{P}\left\{\sum_{k > \gamma_\varepsilon W_\varepsilon} \eta_k^2 s_k^2 > \delta^2\right\} = \mathbf{P}\left\{\sum_{k > \gamma_\varepsilon W_\varepsilon} (\eta_k^2 - 1)s_k^2 > \delta^2 - \sum_{k \geq \gamma W_\varepsilon} s_k^2\right\}$$

$$\leq \mathbf{P}\left\{\sum_{k > \gamma_\varepsilon W_\varepsilon} (\eta_k^2 - 1)s_k^2 > C\varepsilon^2 W_\varepsilon \gamma_\varepsilon^{-2\beta+2}\right\}$$

$$\leq \exp(-C\gamma_\varepsilon W_\varepsilon). \qquad \square$$

A. S. Dalalyan
A. B. Tsybakov
Université Paris VI
2 Place Jussieu, Case 188
75252 Paris Cedex 05
France
E-mail: dalalyan@ccr.jussieu.fr
        tsybakov@ccr.jussieu.fr

G. K. Golubev
Université Aix-Marseille 1
39 rue F. Joliot-Curie
13453 Marseille
France
E-mail: golubev@gyptis.univ-mrs.fr